\documentclass[11pt,reqno]{amsart}

\usepackage[margin=0.5in]{geometry}

\usepackage{amssymb,cite}
\catcode`\@=11 \@addtoreset{equation}{section}

\catcode`\@=12

\usepackage{latexsym}
\usepackage{textcomp}
\usepackage{amsmath}
\usepackage{amsfonts}
\usepackage{amssymb}
\usepackage{mathrsfs}
\usepackage{enumerate}
\usepackage{hyperref}
\renewcommand{\a }{\alpha }
\renewcommand{\b }{\beta }
\renewcommand{\d}{\delta }

\newcommand{\e }{\varepsilon }
\newcommand{\g }{\gamma}

\newcommand{\n }{\nabla }
\newcommand{\vp }{\varphi }
\renewcommand{\phi}{\varphi}

\renewcommand{\o }{\omega }

\newcommand{\be}{\begin{equation}}
\newcommand{\ee}{\end{equation}}

\newcommand{\R}{\mathbb{R}}

\newcommand{\de}{\partial}

\newcommand{\ti}{\widetilde}

\newcommand{\ra}{{\rangle}}
\newcommand{\la}{{\langle}}

\newcommand{\calN }{\mathcal{N}}

\newcommand{\calD }{\mathcal{D}}

\newcommand{\calF}{{\mathcal F}}
\newtheorem{Theorem}{Theorem}[section]

\newtheorem{Lemma}[Theorem]{Lemma}

\newtheorem{Corollary}[Theorem]{Corollary}
\newtheorem{Remark}[Theorem]{Remark}

\def\proof{\noindent{{\bf Proof. }}}
\def\square{\vbox{
    \hrule height .4pt
    \hbox{\vrule width .4pt height 7pt \kern 7pt
       \vrule width .4pt}
    \hrule height .4pt }}

\def\QED{\hfill {$\square$}\goodbreak \medskip}

\def\R{{\mathbb R}}

\def\eps{{\varepsilon}}

\newcommand{\Ds}{(-\Delta)^s}

\newcommand{\dive }{\mathop{\rm div}}

\title[Fractional Schr\"odinger ]{Ground states and concentration
phenomena \\ for the fractional Schr\"odinger equation}

\author{Mouhamed Moustapha Fall}
\address{\noindent M. M. Fall - African Institute for Mathematical Sciences of
Senegal,
AIMS-Senegal
 KM 2, Route de Joal,
 BP:1418,
 Mbour, Senegal.} \email{mouhamed.m.fall@aims-senegal.org}

\author{Fethi Mahmoudi}
\address{\noindent F. Mahmoudi - Departamento de  Ingenieria Matem\'atica
and CMM, Universidad de Chile, Casilla 170 Correo 3, Santiago,
Chile.}\email{fmahmoudi@dim.uchile.cl}

\author{Enrico Valdinoci}
\address{\noindent E. Valdinoci - Weierstra{\ss}
Institut f\"ur Angewandte Analysis und Stochastik,
Mohrenstra{\ss}e 39, 10117 Berlin, Germany.}
\address{\noindent E. Valdinoci -
Dipartimento di Matematica Federigo Enriques,
Universit\`a degli Studi di Milano,
Via Saldini 50, I-20133 Milano, Italy.}
\address{\noindent E. Valdinoci -
Istituto di Matematica Applicata e Tecnologie Informatiche
Enrico Magenes,
Consiglio Nazionale delle Ricerche
Via Ferrata 1, I-27100 Pavia, Italy.}
\email{enrico@math.utexas.edu}

\begin{document}

\date{}
\maketitle

\bigskip

\noindent {\footnotesize{\bf Abstract.}  \sf 
We consider here solutions of the
nonlinear fractional Schr\"odinger equation $$\e^{2s}\Ds u+V( x)u=u^p.$$
We show that concentration points must be critical points for $V$. We  also prove that, if the potential $V$ is coercive and has a unique global minimum,
then ground states concentrate suitably at such minimal point as $\e$ tends to zero.
In addition, if the potential $V$ is radial and radially decreasing,
then the minimizer is unique provided $\e$ is small.
}

\bigskip

\


\

\section{Introduction}\label{s:i}

In this paper we will study standing waves for a nonlinear differential equation driven by the fractional Laplacian. We will focus on the so-called \emph{fractional Schr\"{o}dinger equation}
\begin{equation}\label{eq:1}
\mathrm{i} \hbar\,\frac{\partial \psi}{\partial t} =
\hbar^{2s}\Ds \psi + V(x) \psi - \left|\psi \right|^{p-1} \psi
\end{equation}
where $\hbar$ is the
Planck constant, $(x,t) \in \mathbb{R}^N
\times (0,+\infty)$, $0<s<1$, and $V \colon \mathbb{R}^N \to \mathbb{R}$ is an external potential function. The operator $\Ds$ is the fractional Laplacian of order $s$, which,
for a function
$\vp\in C^\infty_c$ (here and in the sequel when omitting the space of definition  we are meaning $\R^N$) may be defined via
Fourier transform:
$$
\calF{\Ds \vp}(\xi)=|\xi|^{2s} \widehat \vp (\xi) \qquad \text{for
  $\xi \in \R^N$,}
$$
where we used the standard notation
$$
\widehat{\vp}(\xi):=\calF(\vp)(\xi):=\frac{1}{(2\pi)^{\frac{N}{2}}}\,
\int_{\R^N}e^{-\imath\xi\cdot{x} }\vp(x)dx
$$
for the Fourier transform of a function~$\vp\in L^2$.
As customary, we will focus on the standing wave situation
of equation \eqref{eq:1}, namely
on the case in which~$\psi(x,t)=u(x) e^{\frac{\mathrm{i} t}{\hbar}}$, with~$u\ge0$:
under this further assumption (and replacing~$V+1$
with~$V$ and $\hbar$ with the small parameter~$\eps>0$), equation~\eqref{eq:1}
reduces to
\be\label{MAIN:EQ}
\e^{2s}\Ds u+V( x)u-u^p=0.\ee
This is the main equation studied in this
paper and it will be set in the whole of~$\R^N$,
with~$N>2s$ and~$p$ subcritical\footnote{When~$N\le 2s$, one
can say that~$p$ is subcritical when~$p\in(1,+\infty)$.}, namely
\be\label{P HYP}
1<p<\frac{N+2s}{N-2s}.
\ee
As for the potential~$V$ in~\eqref{MAIN:EQ}, we suppose that is smooth, positive,
and bounded from zero, namely we assume that
\be\label{eq:Ass-V}
\|V\|_{C^2}<\infty,\qquad \bar{V}= \inf_{\R^N}V>0. \ee
The weak formulation of the fractional Laplacian
naturally leads to the study of the fractional
Sobolev spaces \be H^s:=\left\{ u\in
L^2\,:\,\int_{\R^N}|\xi|^{2s} |\widehat{u}|^2\,d\xi<\infty  \right\},
\ee endowed with the norm
\begin{eqnarray*}
&& \|u\|_{H^s}^2:=\|u\|_{L^2}^2 +\|u\|_{\calD^{s,2}}^2,\\
{\mbox{where }}&&\|u\|_{\calD^{s,2}}^2:=\int_{\R^N} |\xi|^{2s} |\widehat{u}|^2\,d\xi.\end{eqnarray*}
Notice that all the functional spaces~$L^2$, $H^s$ etc. are set in the
whole of~$\R^n$ unless explicitly mentioned.
In this functional setting, a weak solution
of equation~\eqref{MAIN:EQ} is a function~$u_\e\in H^s$
such that
$$
\e^{2s}\int_{\R^N} |\xi|^{2s} \widehat{u}_\e(\xi)
\overline{\widehat{\varphi}(\xi)}\,d\xi
=\int_{\R^N}
\Big(V( x)u_\e(x)-u_\e^p(x)\Big)\,\varphi(x)\,dx$$
for any~$\varphi\in H^s$. For the  existence
of weak solutions for special cases of~\eqref{MAIN:EQ},
see e.g.~\cite{felmer, dipierro, secchi,CZ,GF,Cho,Cheng,Feng}: in this circumstance,
the solutions found are indeed positive, bounded and $C^{2,\alpha}$
(see Theorem~3.4 in~\cite{felmer}
and Lemma~4.4 in~\cite{cabre-sire}).
In this case,
equation~\eqref{MAIN:EQ}
holds pointwise and the fractional Laplace of~$u$
at the point~$x\in\R^N$ has the integral representation
\be\label{int:re}
(-\Delta)^s u(x)=c(N,s)\int_{\R^N} \frac{2u(x)-u(x+y)-u(x-y)}{|y|^{N+2s}}\,dy
\ee
for a suitable~$c(N,s)>0$, see e.g.
Proposition~3.3 in~\cite{guide}.

The first result that we provide characterizes
the points at which solutions of~\eqref{MAIN:EQ}
concentrate for small~$\eps$, stating that these points
are critical for the potential.  This is somehow an extension to the nonlocal setting\footnote{As a technical comment,
we point out that the proof of Lemma~4.2 in~\cite{Wang}
uses the nondegeneracy of the limit profile, which does not
seem to be available in the fractional setting.
Nevertheless, this nondegeneracy plays a crucial role
only in proving the uniqueness of the maximal points
of the spikes, so we will be able to get around this argument
in our framework.}
of Wang's result, see \cite{Wang}.  
To state this first result,
given a sequence of positive solutions~$u_\e$
for equation~\eqref{MAIN:EQ} in the whole of~$\R^N$, we say that~$x_0\in\R^N$
is a strong concentration point for this sequence (or that
the sequence~$u_\e$ strongly concentrates at~$x_0$) if
\begin{equation}\label{eq:concent}
\begin{gathered}
{\mbox{for any
$\d>0$ there exist~$\e_0$ and~$R>0$ such that, for any $\e\in(0,\e_0)$,}}
\\ {\mbox{$u_\e(x)\leq \d$ for all
$x\in \R^N\setminus B(x_0,\e R)$.}}\end{gathered}
\end{equation}
With this setting, the following result holds:

\begin{Theorem}\label{TH:CON}
Assume \eqref{eq:Ass-V} and let $u_\e \in H^s$ be a sequence of
positive solutions of~\eqref{MAIN:EQ}
in the whole of~$\R^N$ that strongly
concentrate at $x_0$. Then $\n V(x_0)=0$.
\end{Theorem}

We remark that, if we perform a translation and a 
spacial dilation of factor $1/\e$,
equation \eqref{MAIN:EQ} becomes
\be\label{MAIN:EQ:resc}
\Ds u+V( \e x+x_0)u-u^p=0.\ee
Thus, to study the concentration phenomena of
this equation, it is convenient
to define
\begin{equation*}\begin{split}
& V_\e(x):=V( \e x+x_0),\\
& \|u\|_{\e,V}^2:=\int_{\R^N} V_\e(x)\,u^2(x)\,dx,\\
& \|u \|_\e^2:=\|u\|_{\calD^{s,2}}^2+ \|u\|_{\e,V}^2,\\
{\mbox{and }}\quad& \nu(V_\e):=\inf_{{u\ne0}}\frac{\|u\|^2_\e}{
\|u\|_{L^{p+1}}^{2}}.
\end{split}\end{equation*}
Notice that these definitions
also make sense when $\e=0$,
namely, one has
$$  \|u\|_{0,V}^2:=V(x_0)\, \int_{\R^N} u^2(x)\,dx $$
and so on.
Moreover, we remark that if~$u$ is a minimizer for~$\nu(V_\e)$
then~$u_\e(x):=u((x-x_0)/\e)$ is a minimizer for
$$ \nu_\e(V):=\e^{\frac{N\,(1-p)}{1+p}}\inf_{u\ne0} 
\frac{\e^{2s}\|u\|_{\calD^{s,2}}^2+\displaystyle\int_{\R^N}V(x)\,u^2(x)\,dx
}{\|u\|_{L^{p+1}}^2}.$$
In this setting, we can
better determine the variational properties of the concentration point~$x_0$.
Namely, while we know from Theorem~\ref{TH:CON}
that~$x_0$ is a stationary point for the potential, now we give conditions
under which it is a minimum. For this scope,
given a sequence of positive solutions~$u_\e$
for equation~\eqref{MAIN:EQ} in the whole of~$\R^N$, we say that~$x_0\in\R^N$
is a weak concentration point for this sequence (or that
the sequence~$u_\e$ weakly concentrates at~$x_0$) if
there exists a sequence of points~$x_\e\to x_0$ such that
\begin{equation}\label{eq:concent:weak}
\begin{gathered}
{\mbox{for any
$\d>0$ there exist~$\e_0$ and~$R>0$ such that, for any $\e\in(0,\e_0)$,}}
\\ {\mbox{$u_\e(x)\leq \d$ for all
$x\in \R^N\setminus B(x_\e,\e R)$.}}\end{gathered}
\end{equation}
By comparing~\eqref{eq:concent}
and~\eqref{eq:concent:weak}, we notice that strong concentration
implies weak concentration (by choosing~$x_\e:=x_0$) for every $\e$.
Then, the following result holds:

\begin{Theorem}\label{TH:CON:2} 
Suppose that $V$ has a unique global 
minimum point and that $u_\e$ is a minimizer for $\nu_\e(V)$. 
Assume in 
addition that~$V$ at infinity stays above such minimal value, i.e. 
\begin{equation}\label{GHJ} 
\liminf_{|x|\to+\infty} V(x) > \min_{\R^N} 
V.\end{equation} 
Then $u_\e$ weakly concentrates at the global minimum point~$x_0$ of $V$.
More precisely, 
the point~$x_\e$ in~\eqref{eq:concent:weak} is the
unique global maximum point of $u_\e$. \end{Theorem}

We emphasize that, in the above theorem,
an additional complication is that 
the nonlocal operator $\Ds$ does not ``see'' local maximum points.
Namely if $y_\e$ 
is  a local maximum point for   $u_\e$,  it is not necessarily true
that $\Ds u_\e(y_\e)\geq 0$ (and, as a matter
of fact, the ``local'' behavior
of ``nonlocal'' equations can be very wild: for instance
all functions are locally $s$-harmonic up to an arbitrarily small
error, see~\cite{sDSV}). This feature makes the proof
of the uniqueness of the global maximum point of $u_\e$ more delicate than
in the classical case. About characterization of concentration sets for minimizers of singular perturbation problems we refer the reader to 
\cite{AM,Wang,NiTa,NiTa93,SW,NW,DdW} and some references therein.  

Next results establish a uniqueness property for the minimizers:

\begin{Theorem}\label{th:uniq-rad}
Assume that $V\in C^1(\R^N)$, 
with $\inf\limits_{\R^N} V>0$ and it is radial. 
Let~$v_\e$ be a minimizer for $$
\frac{\e^{2s}\|u\|_{\calD^{s,2}}^2+\displaystyle\int_{\R^N}V(x)\,u^2(x)\,dx
}{\|u\|_{L^{p+1}}^2}$$
in the class of radial competitors~$u\in H^s$, $u\ne0$.
Then $v_\e$
is unique, provided that~$\e$ small enough.
\end{Theorem} 

\begin{Corollary}\label{coroll:uniq-rad}
Assume that $V\in C^1(\R^N)$,
with $\inf\limits_{\R^N} V>0$ and it is radial
and radially decreasing. 
Let~$v_\e$ be a minimizer for $\nu_\e(V)$.
Then $v_\e$
is unique, provided that~$\e$ small enough.
\end{Corollary}

The rest of the paper is organized as follows. 
In Section~\ref{SSS1}
we study the concentration phenomena at given points of the
space and we prove Theorem~\ref{TH:CON}. The proof of
Theorem~\ref{TH:CON:2} requires some preliminary work,
that is carried out in Section~\ref{SSS8}. In particular,
we obtain there an expansion of the minimizers of~$\nu(V_\e)$
as perturbation of a suitable translation of the ground state
(for this, no condition on the concentration point is required).

The proof of Theorem~\ref{TH:CON:2} is then completed in Section~\ref{CCC0}.
Then, Section~\ref{PR1} contains the preliminaries needed for
the proof of Theorem~\ref{th:uniq-rad}, which, in turn, will
be completed in Section~\ref{PR2}. The proof of
Corollary \ref{coroll:uniq-rad}
is contained in Section~\ref{PR2bis}.
\bigskip

\subsection*{Acknowledgments}
It is a pleasure to thank
Rupert Frank for his very deep comments on an earlier
version of this manuscript.

MMF is supported by the Alexander von Humboldt Foundation
and partially by the Simons Associateship funding from  the International Center for Theoretical Physics (ICTP). 
FM is supported by the Fondecyt grant 1140311, Fondo Basal CMM and  partially by ``Millennium Nucleus Center for Analysis of PDE NC130017".
EV is supported by 
the ERC grant $\e$ ``Elliptic
Pde's and Symmetry of Interfaces and Layers for Odd Nonlinearities'' and
the PRIN grant 
``Aspetti variazionali e
perturbativi nei problemi differenziali nonlineari''.

\section{Concentrations occurring at critical points of $V$
and proof of Theorem~\ref{TH:CON}}\label{SSS1}
In this section, we prove Theorem~\ref{TH:CON}.
We define
$$
v_\e(x):=u_\e(\e x+x_0)
.$$
By construction, $v_\e$ is a positive solution of
\begin{equation}\label{eq:vesatisf}
\Ds v_\e+V(\e x+x_0)v_\e-v_\e^p=0\qquad {\mbox{in }}\R^N.
\end{equation}
Roughly speaking, the idea is to take the derivative
of~\eqref{eq:vesatisf},  test it against~$v_\e$,
 integrate by parts and hence  send~$\e\to0$, in order to
see that~$\nabla V(x_0)=0$: but to do these steps,
some uniform regularity and decay estimates in~$\e$ are in order.
To obtain these estimates,
we define
$$
m_\e:=\max_{\R^N} v_\e= \|u_\e\|_{L^\infty}.
$$
We claim that
\begin{equation}\label{m-eps}
m:=\sup_{\e\in(0,1)}m_\e<+\infty.
\end{equation}
The proof is based on a classical contradiction
and scaling arguments.
Namely, suppose that
\begin{equation}\label{C1}
m_\e\to+\infty, \end{equation}
up to a subsequence.
Now we recall~\eqref{eq:Ass-V}
and we use~\eqref{eq:concent} with
$$\d:=\min\left\{1,\, \left(\frac{\bar{V}}{2}\right)^{\frac{1}{p-1}}\right\}.$$
Accordingly, we obtain that there exists~$R_1>0$
for which
\be\label{C2}
u_\e(x)\le \min\left\{1,\, \left(\frac{\bar{V}}{2}\right)^{\frac{1}{p-1}}\right\} \
{\mbox{ for any~$y\in\R^N$
such that~$|y-x_0|\ge\e R_1$,}}\ee
as long as~$\e$ is small
enough. Now we notice that if~$x\in \R^N\setminus B(0,R_1)$
then~$|(\e x+x_0)-x_0|\ge \e R_1$: hence~\eqref{C2}
implies that
\be\label{C3}
v_\e(x)\le \min\left\{1,\, \left(\frac{\bar{V}}{2}\right)^{\frac{1}{p-1}}\right\} \
{\mbox{ for any~$x\in\R^N\setminus B(0,R_1)$.}}\ee
{F}rom~\eqref{C1} and~\eqref{C3}, we conclude that, for small~$\e$,
$$ 1<m_\e=\max_{\overline{B(0,R_1)}}v_\e,$$
and there exists
\be\label{9u}
x_\e\in \overline{B(0,R_1)}\ee
maximizing~$v_\e$, that is
$$ m_\e = \max_{\R^N} v_\e = v_\e(x_\e).$$
So, we set $\mu_\e:=m_\e^{\frac{1-p}{2s}}$
and~$w_\e(x):= m_\e^{-1} v_\e (x_\e+\e \mu_\e x)$. Then~$\|w_\e\|_{L^\infty}=1=w_\e(0)$ and
\begin{equation}\label{Ca0}\begin{split}
\Ds w_\e(x) \,=& -\mu_\e^{2s}  V(\e (x_\e+\e\mu_\e))\,w_\e(x)+
w_\e^p(x).\end{split}\end{equation} Notice that~$\mu_\e\to0$ as~$\e\to0$, thanks to~\eqref{C1}.
Therefore, by~\eqref{Ca0},
we have that~$\| \Ds w_\e\|_{L^\infty}$ is bounded uniformly in~$\e$. As a
consequence of this and of the regularity results (see e.g.
Lemma~4.4 in~\cite{cabre-sire}, see also \cite{Sil}), we deduce
that~$\|w_\e\|_{C^{2,\alpha}}$ is bounded uniformly in~$\e$, for
some~$\alpha\in(0,1)$. Hence, we can suppose that~$w_\e$ converges
to some function~$w_0$ in~$C^{2,\alpha}_{loc}$, with
\be\label{77}
\|w_0\|_{L^\infty}=1=w_0(0).\ee
By passing to the limit in~\eqref{Ca0}, we obtain
that
\be\label{777}
\Ds w_0= w_0^p \ {\mbox{ in $\R^N$.}}\ee
Since the only non-negative and bounded solution
of~\eqref{777} with~$p$ subcritical
(according to~\eqref{P HYP})
is the one constantly equal to zero (see Remark~1.2
in~\cite{JLX} or Theorem~1.3
in~\cite{FW}),
we conclude that~$w_0$ vanishes identically,
in contradiction with~\eqref{77}.

This completes the proof of~\eqref{m-eps}.
As a consequence of~\eqref{eq:vesatisf}, \eqref{m-eps}
and of the regularity results
(see e.g. Lemma~4.4 in~\cite{cabre-sire}), we conclude that
\be\label{c 2 a 0}
{\mbox{$\|v_\e\|_{C^{2,\alpha}}$
is bounded uniformly in~$\e$}},\ee
hence we may suppose that
\be\label{3.9bis}
{\mbox{$v_\e$ converges to some function~$v_0$ in~$C^{2,\alpha}_{loc}$.}}\ee
Now, since~$x_\e$ maximizes~$v_\e$, we have that
$$ 2v_\e(x_\e)-v_\e(x_\e+y)-v_\e(x_\e-y)\ge0 \quad
{\mbox{ for any }} y\in\R^N,$$
and so, using \eqref{int:re} and~\eqref{eq:vesatisf},
$$ 0\le (-\Delta)^s v_\e(x_\e)= v_\e^p (x_\e)-
V(\e x_\e+x_0)v_\e(x_\e).$$
Accordingly,
\be\label{C9}
V(\e x_\e+x_0)\le v_\e^{p-1} (x_\e).
\ee
Since~$|x_\e|$ is bounded uniformly in~$\e$, in light of~\eqref{9u},
we suppose, up to a subsequence, that~$x_\e\to \bar x$,
for some~$\bar x\in \overline{B(0,R_1)}$, as~$\e\to0$.
Thus, by taking the limit as~$\e\to0$
in~\eqref{C9},
we obtain that
$$ 0<\bar{V}\le V(x_0)\le v_0^{p-1} (\bar x).$$
In particular,
\be\label{n i z}
{\mbox{$v_0$ is not identically zero.}}
\ee
Next we claim that there exists $\e_0>0$ such that
\be\label{po de}
v_\e(x)\leq \frac{Const}{1+|x|^{N+2s}}\qquad \forall \e\in (0,\e_0).
\ee
To prove this,
we use Lemma~4.2 of~\cite{felmer}, according to which there
exists a function~$\bar w$ such that
\be\label{Qc} 0\le\bar w(x)\le \frac{Const}{1+|x|^{N+2s}}\ee
and
\be\label{Cq} (-\Delta)^s \bar w+\frac{\bar{V}}{2} w\ge 0 \
{\mbox{ in }}\R^N\setminus B(0,\tilde R),\ee
for a suitable~$\tilde R>0$.
Now, we take
\be\label{eb}
R_2:=\min\{R_1,\,\tilde R\},\ee
where~$R_1$ is the one in~\eqref{C3}.
Thanks to \eqref{eq:concent}, we have that $v_\e$ converges to zero as $|x|\to \infty$ uniformly with respect to $\e$. {F}rom~\eqref{C3}, we obtain
\begin{equation}\begin{gathered}\label{86}
(-\Delta)^s v_\e +\frac{\bar{V}}{2} v_\e
=(-\Delta)^s v_\e + V v_\e
-\left(V-\frac{\bar{V}}{2}\right) v_\e \\
=v_\e^p
-\left(V-\frac{\bar{V}}{2}\right) v_\e\le
v_\e^p-\frac{\bar{V}}{2} v_\e= v_\e \left( v_\e^{p-1}-\frac{\bar{V}}{2}\right)\le 0.
\end{gathered}\end{equation}
Now we set
\begin{equation}\label{2.18BIS}
b:=\inf_{B(0,R_2)} \bar w>0\end{equation}
and
\be\label{091}
z_\e:=(m+1) \bar w-bv_\e,\ee where~$m$ is given in~\eqref{m-eps}.
Our goal is to show that
\be\label{090}
z_\e \ge0 \ {\mbox{ in }} \R^N.\ee
For this we argue by contradiction and suppose that
\be\label{QQ}
0>\inf_{\R^N} z_\e =\lim_{j\to+\infty} z_\e (x_{j,\e}),\ee
for a suitable sequence~$x_{j,\e}$. Notice that
$$ \lim_{|x|\to+\infty} \bar w(x)=0,$$
due to~\eqref{Qc}, and
$$ \lim_{|x|\to+\infty} u_\e(x)=0,$$
due to our integrability and continuity assumptions on~$u_\e$,
and therefore
$$ \lim_{|x|\to+\infty} v_\e(x)=0,$$
and so
$$ \lim_{|x|\to+\infty} z_\e(x)=0.$$
Consequently, the sequence~$x_{j,\e}$
is bounded and therefore, up to subsequence, we suppose
that~$x_{j,\e}\to x_{\star,\e}$ as~$j\to+\infty$,
for some~$x_{\star,\e}\in\R^N$. So~\eqref{QQ} becomes
\be\label{QQ2} 0>\min_{\R^N} z_\e =z_\e (x_{\star,\e}).\ee
The minimality property of~$x_{\star,\e}$ and~\eqref{int:re}
give that
\be\label{QQ3}
(-\Delta)^s z_\e(x_{\star,\e})=
c(N,s)\int_{\R^N} \frac{2z_\e(x_{\star,\e})-
z_\e(x_{\star,\e}+y)-z_\e(x_{\star,\e}-y)}{|y|^{n+2s}}\,dy\le 0
.\ee
Now notice that, by \eqref{m-eps} and~\eqref{2.18BIS},
$$ z_\e \ge mb + \bar w-bm>0 \ {\mbox{ in }} B(0,R_2).$$
Comparing this with~\eqref{QQ2}, we see that
\be\label{PP} x_{\star,\e}\in \R^N\setminus B(0,R_2).\ee
Moreover, from~\eqref{Cq}, \eqref{eb} and~\eqref{86},
we obtain that
\be\label{092} (-\Delta)^s z_\e +\frac{\bar{V}}{2} z_\e \ge0 \
{\mbox{ in }} \R^N\setminus B(0,R_2).\ee
Thanks to~\eqref{PP}, we can evaluate~\eqref{092} at the point~$x_{\star,\e}$:
in this way, and recalling~\eqref{QQ2} and~\eqref{QQ3},
we obtain that
$$ 0\le (-\Delta)^s z_\e(x_{\star,\e}) +\frac{\bar{V}}{2} z_\e(x_{\star,\e})
<0.$$
This is a contradiction, so~\eqref{090} is established. 

{F}rom~\eqref{090},
we deduce that~$v_\e\le (m+1) b^{-1} \bar w$, which, together
with~\eqref{Qc}, completes the proof of~\eqref{po de}.

Using~\eqref{3.9bis} and~\eqref{po de} and the
dominated convergence theorem, we see that
\be\label{L2.1}
\lim_{\e\to0}\int_{\R^N} \partial_i V(\e x+x_0) v_\e^2 =\partial_i V(x_0)
\int_{\R^N}v_0^2,
\ee
for any~$i\in\{1,\dots,N\}$.

Now we show that
\begin{equation}\label{ADD}
|\nabla v_\e| \in L^2.
\end{equation}
For this, we write \eqref{eq:vesatisf} as
\begin{equation}\label{wfd34efd}
(-\Delta)^s v_\e = -V(\e x +x_0)\,v_\e+v_\e^p=:\psi_\e.\end{equation}
We know from \eqref{po de} that
$$ \psi_\e(x)\le\frac{Const}{1+|x|^{N+2s}}.$$
Thus we take the fundamental solution~$\Gamma$ of the operator
in~\eqref{wfd34efd} and we obtain that
$$ v_\eps(x)=\,Const\,\int_{\R^N} \frac{\psi_\e(y)}{|x-y|^{N-2s}}\,dy.$$
Therefore
$$ |\nabla v_\eps(x)|\le Const\,
\int_{\R^N} \frac{|\psi_\e(y)|}{|x-y|^{N-2s+1}}\,dy\le
\int_{\R^N} \frac{Const}{(1+|y|^{N+2s})\,|x-y|^{N-2s+1}}\,dy.$$
So, fixing~$x\in\R^N\setminus B(0,2)$, we observe that
$$ \int_{B(x,1)} \frac{|\psi_\e(y)|}{|x-y|^{N-2s+1}}\,dy\le
\frac{Const}{|x|^{N+2s}}
\int_{B(0,1)} \frac{1}{|\xi|^{N-2s+1}}\,d\xi
\le\frac{Const}{|x|^{N+2s}},$$
therefore we obtain that
\begin{eqnarray*}
|\nabla v_\eps(x)|&\le& Const\,\left[
\int_{\R^N\setminus B(x,1)} \frac{1}{(1+|y|^{N+2s})\,|x-y|^{N-2s+1}}\,dy
+\frac{1}{|x|^{N+2s}}\right]\\
&\le&
Const\,\left[
\int_{\R^N\setminus B(x,1)} \frac{1}{(1+|y|^{N+2s})\,(1+|x-y|^{N-2s+1})}\,dy
+\frac{1}{|x|^{N+2s}}\right]
.\end{eqnarray*}
Accordingly, by the properties of the convolution of decaying kernels
(see e.g. Lemma~5.1 in~\cite{davila-delpino-dipierro-valdinoci}), we obtain that
\begin{equation}\label{6-001}
|\nabla v_\eps(x)|\le\frac{Const}{|x|^{\kappa}},
\end{equation}
with~$\kappa:=\min\{ N+2s,\; N-2s+1\}$.
Notice that
$$ 2\kappa=\min\{ 2N+4s,\; N+N-4s+2\}>
\min\{2N,\; 2s+N-4s+2\}>N,$$
hence~\eqref{6-001}
implies~\eqref{ADD}, as desired.

Now we perform some calculations on integrals
that involve~$v_\e$.
For this, we let~$e_i$ be the $i$th vector
of the standard Euclidean base, we fix~$R>1$
and we use the divergence theorem to see that,
for any~$i \in\{1,\dots,N\}$,
$$ \int_{B(0,R)} \partial_i v_\e^{p+1}=
\int_{B(0,R)} \dive ( v_\e^{p+1} e_i)=
\int_{\partial B(0,R)} v_\e^{p+1} \frac{x_i}{R}.$$
Thus, from~\eqref{po de}, we have
\be\label{I E 1}
\int_{B(0,R)} \partial_i v_\e^{p+1}=O(R^{N-1-(p+1)(N+2s)}).\ee
Similarly,
\begin{eqnarray*} && \int_{B(0,R)}
V(\e x+x_0) \partial_i v_\e^2=
\int_{B(0,R)} \Big[ \dive \Big( V(\e x+x_0) v_\e^2 e_i\Big)
- \e \partial_i V(\e x+x_0) v_\e^2 \Big]\\&&\qquad
= \int_{\partial B(0,R)} V(\e x+x_0) v_\e^2 \frac{x_i}{R}
- \e \int_{B(0,R)} \partial_i V(\e x+x_0) v_\e^2
\end{eqnarray*}
which, together with~\eqref{po de}, gives that
\be\label{I E 2}
\int_{B(0,R)}
V(\e x+x_0) \partial_i v_\e^2 =
- \e \int_{B(0,R)} \partial_i V(\e x+x_0) v_\e^2
+O(R^{N-1-2(N+2s)}).\ee
We summarize the estimates in~\eqref{I E 1} and~\eqref{I E 2}
by writing
\be\label{IE}\begin{split}
&\int_{B(0,R)}
\Big( \e\partial_i V(\e x+x_0)v_\e^2
+\frac12 V(\e x+x_0) \partial_i v_\e^2 -\frac{p}{p+1}\partial_i v_\e^{p+1}\Big)
\\ &\qquad= \frac\e2 \int_{B(0,R)}\partial_i V(\e x+x_0)v_\e^2+O(R^{N-1-2(N+2s)}).
\end{split}\ee
Now, we point out that~$(-\Delta)^s v_\e$ is $C^{2}$
and bounded,
due to~\eqref{eq:Ass-V},
\eqref{eq:vesatisf}, \eqref{m-eps} and~\eqref{c 2 a 0}
(recall also~\eqref{eq:Ass-V}, therefore
we can speak
about~$\partial_i (-\Delta)^s v_\e$ in the classical sense.
Accordingly, we can take a derivative, say in the $i$th coordinate
direction, of~\eqref{eq:vesatisf}: we get
\be\label{d i}
\partial_i \Ds v_\e+
\e\partial_i V(\e x+x_0)v_\e
+V(\e x+x_0)\partial_i v_\e -p v_\e^{p-1} \partial_i v_\e=0.
\ee
So, recalling~\eqref{c 2 a 0}, \eqref{po de} and~\eqref{ADD}, we see that
\be\label{P00}
\partial_i \Ds v_\e\in L^2.\ee
Consequently, by Plancherel theorem, we obtain
\begin{equation}\label{P0}
\begin{split}
&\int_{\R^N} v_\e \partial_i \Ds v_\e
= \int_{\R^N}\overline{\widehat{v_\e}}\,\calF({\partial_i \Ds v_\e})\\
&\qquad= -\int_{\R^N}\xi_i \overline{\widehat{v_\e}}\,\calF({\Ds v_\e})= -\int_{\R^N} \xi_i |\xi|^{2s}|\widehat{v_\e}(\xi)|^2\\
&\qquad= -\int_{\R^N} \xi_i \big| \calF \big( {(-\Delta)^{s/2} v_\e }\big)\big|^2
.\end{split}\end{equation}
We remark that
\be\label{98}
(-\Delta)^{s/2} v_\e\in L^2.\ee
Indeed, since~$u_\e\in H^s$ (hence~$v_\e\in H^s$),
we have that
\begin{eqnarray*}
\int_{\R^N} |\xi|^{2s} |\widehat{v_\e}|^2\le
\int_{\R^N} (1+|\xi|^{2s}) |\widehat{v_\e}|^2<+\infty,
\end{eqnarray*}
therefore~$\calF({(-\Delta)^{s/2} v_\e})=|\xi|^s \widehat{v_\e}\in L^2$
and so~\eqref{98} follows from the 
Plancherel theorem.

Also, for any~$g\in L^2$, we have that
\be\label{98.b}
{\mbox{the map $\xi\mapsto|\widehat{g}(\xi)|^2$ is even.}}
\ee
To check this, notice that
$$ \overline{\widehat{g}(\xi)}=
\overline{\int_{\R^N} g(x) e^{-\mathrm{i} x\cdot\xi} \,dx}
=\int_{\R^N} g(x) e^{\mathrm{i} x\cdot\xi} \,dx = \widehat{g}(-\xi)$$
and so
$$ \overline{\widehat{g}(-\xi)}=
\widehat{g}(\xi).$$
As a consequence
$$ |\widehat{g}(-\xi)|^2=
\overline{\widehat{g}(-\xi)} \widehat{g}(-\xi)=
\widehat{g}(\xi) \overline{\widehat{g}(\xi)}=|\widehat{g}(\xi)|^2,$$
that proves~\eqref{98.b}.

So, from~\eqref{98} and~\eqref{98.b}, we see
that the map~$\xi\mapsto \xi_i \,| \calF({(-\Delta)^{s/2} v_\e })|^2$
is odd, and therefore
$$ \int_{B(0,R)} \xi_i \,\big| \calF\big({(-\Delta)^{s/2} v_\e }\big)\big|^2=0$$
for any~$R>0$.
By plugging this into~\eqref{P0} and recalling~\eqref{P00} we obtain
\be\label{P1} \int_{\R^N} v_\e \partial_i \Ds v_\e
=\lim_{R\to+\infty}
\int_{B(0,R)} v_\e \partial_i \Ds v_\e=0.\ee
Now we go back to~\eqref{d i} and we multiply
this equation by~$v_\e$: in this way we
obtain that
\begin{eqnarray*}
0&=& v_\e \partial_i \Ds v_\e+
\e\partial_i V(\e x+x_0)v_\e^2
+V(\e x+x_0)v_\e\partial_i v_\e -p v_\e^p \partial_i v_\e
\\ &=&
v_\e \partial_i \Ds v_\e+\e\partial_i V(\e x+x_0)v_\e^2
+\frac12 V(\e x+x_0) \partial_i v_\e^2 -\frac{p}{p+1}\partial_i v_\e^{p+1}.
\end{eqnarray*}
We fix~$R>1$ and we integrate the above
equation on~$B(0,R)$: thus, exploiting~\eqref{IE}
we obtain
$$ \int_{B(0,R)} v_\e \partial_i \Ds v_\e+
\frac{\e}{2} \int_{B(0,R)}\partial_i V(\e x+x_0)v_\e^2 = O(
R^{N-1-2(N+2s)}).$$
So we send~$R\to+\infty$, recalling also~\eqref{P1}
and we divide by~$\e$,  we get
$$ \int_{\R^N}\partial_i V(\e x+x_0)v_\e^2 =0.$$
Now we send~$\e\to0$: recalling~\eqref{L2.1}
we conclude that
$$ \partial_i V(x_0)\int_{\R^N}v_0^2=0.$$
Therefore, by~\eqref{n i z}, we obtain that~$\partial_i V(x_0)=0$
for any~$i\in\{1,\dots,N\}$,
and this completes the proof of Theorem~\ref{TH:CON}.~\QED

\section{Concentration points of ground-states: preliminary work
for the proof of Theorem~\ref{TH:CON:2}}\label{SSS8}
In this section we discuss some basic concentration properties of
the minimizers. For this, we recall that, for any~$\lambda>0$
there exists a unique function~$U_\lambda$ that attains the following minimization
problem
$$ \nu(\lambda):=\inf_{\|u\|_{L^{p+1}}=1} \|u\|_{\calD^{s,2}}^2 +\lambda \|u\|_{L^2}^2.$$
In addition, such minimizer is unique radially symmetric and belongs to $C^\infty\cap H^{2s+1}(\R^N)$ (we refer to~\cite{FLS}
for further details on this, see in particular
Theorem~4 there).
Thus, we will denote by~$\ti{U}$ the radially symmetric function
that attains 
$$ \inf_{\|u\|_{L^{p+1}}=1} \|u\|_{\calD^{s,2}}^2 +\ti{V}\,\|u\|_{L^2}^2,$$
where
\begin{equation}\label{DEFTIV} \ti{V}:=\inf_{\R^N} V.\end{equation}
With this notation, we provide an asymptotic expansion
for the minimizers of~$\nu(V_{\e})$. It is worth pointing out that
this expansion
is valid without assuming any structural condition on the potential~$V$
(in particular the point~$x_0\in\R^n$ can be fixed, without
assuming that is minimal or critical):

\begin{Lemma}\label{no0l}
Let $v_\e$ be a positive
minimizer for $\nu(V_{\e})$, with~$\|v_\e\|_{L^{p+1}}=1$. Then  
there exists a sequence of points $a_\e, c\in\R^N$ 
such that, up to a subsequence,
\begin{eqnarray*}
&& v_\e(x+a_\e)=\ti{U}(x- c)+\o_{\e}(x),\\
{\mbox{with }}&& \|\o_\e\|_{H^{s}}\to 0 \ {\mbox{ as }}\ \e\to0.\end{eqnarray*}
Also
$$\lim_{\e\to0}\nu(V_\e)=\nu(\ti{V})$$
and, for any~$x\in\R^n$,
\begin{equation}\label{DEFTIV-1} 
\lim_{\e\to0}V(\e x+\e a_\e+x_0)
=\ti{V}.\end{equation}
\end{Lemma}

\proof We observe that
$$ \|u\|_{\e,V}^2 \in \Big[ \bar V\,\|u\|_{L^2}^2,\;\|V\|_{L^\infty} \|u\|_{L^2}^2\Big]$$
thanks to~\eqref{eq:Ass-V}, and therefore~$\nu(V_{\e})$ is bounded (and bounded from zero) uniformly in~$\e$.
Hence, up to a subsequence, we suppose that
\begin{equation}\label{LIMI1}
\nu(V_{\e})\to\ti{\nu} \end{equation}
as~$\e\to0$, for some~$\ti{\nu}>0$.

Also, $v_\e$ is bounded in~$H^s$ and,
using Lemma~2.2 in~\cite{felmer}, we have that
there exists $a_\e\in \R^N$ and positive real numbers $R$ and $\g$   such that 
\be \label{eq:compact_infty}
\liminf_{\e\to 0}\int_{B_R(a_\e)}v_\e(x)\,dx\geq \g.
\ee
Thus, setting $w_\e(x)=v_\e(x+a_\e)$, we have
that $w_\e$ is bounded in $H^s$ so
it converges, up to a subsequence, to a function $w\in H^ s$ 
weakly in~$H^s$, strongly in~$L^{p+1}_{loc}$ and a.e.; furthermore,
by \eqref{eq:compact_infty}, we have that~$w\neq 0$. 

We also notice that, by \eqref{eq:Ass-V} and the
theorem of Ascoli,
there exists~$\lambda:\R^N\to\R$ such that,
up to a subsequence, 
\begin{equation}\label{laof01}
\lim_{\e\to0}V(\e x+\e a_\e+x_0)=\lambda(x).
\end{equation}
We set~$\lambda:=\lambda(0)$ and we claim that
\begin{equation}\label{laof0}
\lambda(x)=\lambda
\end{equation}
for any~$x\in\R^N$. Indeed, for any~$x\in\R^N$,
\begin{eqnarray*}
&& |\lambda(x)-\lambda|=\lim_{\e\to0}\big|V(\e x+\e a_\e+x_0)-
V(\e a_\e+x_0)\big|\\
&&\qquad\le \,Const\,\lim_{\e\to0} |\e x| =0,
\end{eqnarray*}
thanks to \eqref{eq:Ass-V}, and this proves~\eqref{laof0}.

By \eqref{laof01} and~\eqref{laof0}, we can write
\begin{equation}\label{dsvcik2ws} \lim_{\e\to0}V(\e x+\e a_\e+x_0)=\lambda.\end{equation}
Since, by \eqref{eq:vesatisf},
$$
\Ds w_\e+ V(\e x+\e a_\e+x_0) w_\e=\nu(V_{\e}) w_\e^{p},
$$
we can pass to the limit and obtain
\be\label{Dsv0}
\Ds w+ \lambda w=\ti{\nu} w^{p}.
\ee
By testing \eqref{Dsv0} against $w$ we obtain that
\be\label{u1}
\ti{\nu}=\frac{
\|w\|_{\calD^{s,2}}^2 +\lambda\,\|w\|_{L^2}^2}{\|w\|^2_{L^{p+1}}} \ge
\inf_{u\ne0}\frac{\|u\|_{\calD^{s,2}}^2 +\lambda\,\|u\|_{L^2}^2}{\|u\|^2_{L^{p+1}}} =\nu(\lambda).
\ee
On the other hand, by the
dominated convergence theorem, we see that, for any~$u\in C^\infty_c$,
$$ \lim_{\e\to0}\int_{\R^N} \big| V(\e x+\e a_\e+x_0)-\lambda\big| \,u^2(x)\,dx=0.$$
As a consequence, for any~$u\in C^\infty_c$, $u\ne0$, we set~$\tilde u_\e(x):=u(x-a_\e)$
and we observe that
\begin{eqnarray*}
\ti{\nu}&=&\lim_{\e\to0} \nu(V_{\e})\\
&\le& \frac{\|u\|_{\calD^{s,2}}^2 + \|u\|_{\e,V}^2}{ \|u\|_{L^{p+1}}^2 }
\\ &=&\lim_{\e\to0}
\frac{\|\tilde u_\e\|_{\calD^{s,2}}^2 + 
\displaystyle\int_{\R^N}V(\e x+x_0)\,\tilde u_\e^2(x)\,dx}{
\|\tilde u_\e\|_{L^{p+1}}^2 }\\
\\ &=&
\lim_{\e\to0}
\frac{\|u\|_{\calD^{s,2}}^2 +
\displaystyle\int_{\R^N}V(\e x+x_0)\,u^2(x-a_\e)\,dx}{
\|u\|_{L^{p+1}}^2 }
\\ &=&
\lim_{\e\to0}
\frac{\|u\|_{\calD^{s,2}}^2 +
\displaystyle\int_{\R^N}V(\e x+\e a_\e+x_0)\, u^2(x)\,dx}{
\|u\|_{L^{p+1}}^2 }
\\ &=&
\frac{\|u\|_{\calD^{s,2}}^2 +\lambda\|u\|_{L^2}^2
}{ \|u\|_{L^{p+1}}^2 } .
\end{eqnarray*}
By density, this is valid for any~$u\in H^s$,
and so, taking the infimum over~$u\ne0$, we obtain that~$\ti{\nu}\le
\nu(\lambda)$. This and~\eqref{u1} give that
\begin{equation}\label{LIMI2}\ti{\nu}=\nu(\lambda).\end{equation}
This, \eqref{Dsv0}
and the uniqueness of the ground state (see Theorem~4 
in~\cite{FLS})
give that $w$ is a translation of $U_\lambda$, namely~$w(x)=U_\lambda(x-c)$,
for some~$c\in\R^N$.

Now we claim that
\begin{equation}\label{LIMI3}
\lambda=\ti{V}
\end{equation}
To prove this, let us fix~$q\in\R^n$. Then, for any~$u\in C^\infty_c$, 
we set~$\underline u_\eps(x):=u(x+\eps^{-1}(x_0-q))$ and we use
the change of variable~$y:=x+\e^{-1}(q-x_0)$ to obtain that
\begin{eqnarray*}
&& \frac{\|u\|_{\calD^{s,2}}^2 +
\displaystyle\int_{\R^N}V(\e x+q)\,u^2(x)\,dx}{\|u\|_{L^{p+1}}^2} \\
&=& \frac{\|u\|_{\calD^{s,2}}^2 +
\displaystyle\int_{\R^N}V(\e y+x_0)\,u^2(y+\e^{-1}(x_0-q))\,dy}{
\|u\|_{L^{p+1}}^2} \\
&=& \frac{\|\underline u_\eps\|_{\calD^{s,2}}^2 +
\displaystyle\int_{\R^N}V(\e y+x_0)\,\underline u_\eps^2(y)\,dy}{
\|\underline u_\eps\|_{L^{p+1}}^2} \\
&=& \frac{\|\underline u_\e\|_{\e,V}^2}{\|\underline u_\eps\|_{L^{p+1}}^2}\\
&\ge& \nu(V_\e).
\end{eqnarray*}
So, by \eqref{LIMI1}, \eqref{LIMI2}
and the dominated convergence theorem, we obtain
\begin{eqnarray*}
&& \frac{\|u\|_{\calD^{s,2}}^2 +V(q)\|u\|_{L^2}^2
}{\|u\|_{L^{p+1}}^2} =\lim_{\e\to0}
\frac{\|u\|_{\calD^{s,2}}^2 +
\displaystyle\int_{\R^N}V(\e x+q)\,u^2(x)\,dx}{\|u\|_{L^{p+1}}^2} \\ &&\quad\ge
\lim_{\e\to0}\nu(V_\e)=
\ti{\nu}=\nu(\lambda)=
\frac{ \|U_\lambda\|_{\calD^{s,2}}^2 +\lambda\|U_\lambda\|_{L^2}^2
}{\|U_\lambda\|_{L^{p+1}}^2}.\end{eqnarray*}
This is valid for any~$u\in C^\infty_c$ and so, by density, also for~$U_\lambda$.
Thus we conclude that
$$ \frac{\|U_\lambda\|_{\calD^{s,2}}^2 +V(q)\|U_\lambda\|_{L^2}^2
}{\|U_\lambda\|_{L^{p+1}}^2}\ge
\frac{\|U_\lambda\|_{\calD^{s,2}}^2 +\lambda\|U_\lambda\|_{L^2}^2
}{\|U_\lambda\|_{L^{p+1}}^2}$$
and therefore
$$ V(q)\ge\lambda.$$
Now, this is valid for any~$q\in\R^N$, thus, recalling~\eqref{DEFTIV},
we obtain that
$$ \ti{V}=\inf_{q\in\R^N} V(q)\ge\lambda.$$
The other inequality follows from~\eqref{dsvcik2ws}, and so the proof
of~\eqref{LIMI3} is complete.

Then, \eqref{LIMI3} and the definition of~$\ti{U}$ give that~$U_\lambda=\ti{U}$.
Accordingly,
$$ v_\e(x+a_\e) =w_\e(x)\to w(x)=U_\lambda(x-c)=\ti{U}(x-c)$$
weakly in~$H^s$, strongly in~$L^{p+1}_{loc}$ and a.e., so
to complete the proof of
Lemma~\ref{no0l}
it only remains to show that the convergence occurs strongly in~$H^s$. To see this, we use the
fact that $w$ is a minimizer for the quotient $\nu(\ti{V})=\nu(\lambda)=\ti{\nu}$, hence
\begin{equation} \label{ocxw45676543edfg} \ti{\nu}=\frac{ 
\|w\|_{\calD^{s,2}}^2 +
\lambda\|w\|_{L^2}^2}{\|w\|_{L^{p+1}}^{2}}.\end{equation}
On the other hand, by testing~\eqref{Dsv0} against~$w$, we obtain that
$$ \|w\|_{\calD^{s,2}}^2 +
\lambda\|w\|_{L^2}^2 =\ti{\nu}\,\|w\|_{L^{p+1}}^{p+1}.$$
By comparing this with~\eqref{ocxw45676543edfg}, we conclude that~$
\|w\|_{L^{p+1}}=1 $. 
Therefore
\begin{eqnarray*}
&& \|w_\e\|_{\calD^{s,2}}^2 +
\lambda\|w_\e\|_{L^2}^2=
\|v_\e\|_{\calD^{s,2}}^2 +
\lambda\|v_\e\|_{L^2}^2\\
&&\quad=
\|v_\e\|_{\calD^{s,2}}^2 +
\lambda\|v_\e\|_{\e,V}^2+\int_{\R^N}\big( \lambda-V(\e x+x_0)\big)\,v_\e(x)\,dx
\\&&\quad=\nu(V_\e)+
\int_{\R^N}\big( \lambda-V(\e x+x_0)\big)\,v_\e(x)\,dx.
\end{eqnarray*}
Moreover, by \eqref{DEFTIV} and~\eqref{LIMI3},
$$ \lambda=\inf_{\R^N}V\le V(\e x+x_0),$$
thus we obtain that
$$ \|w_\e\|_{\calD^{s,2}}^2 +
\lambda\|w_\e\|_{L^2}^2 \le \nu(V_\e).$$
So, from the weak convergence and Fatou lemma,
passing to the limit we obtain that
\begin{eqnarray*}
&&\ti{\nu}=\|w\|_{\calD^{s,2}}^2 +
\lambda\|w\|_{L^2}^2\le 
\liminf_{\e\to0}
\|w_\e\|_{\calD^{s,2}}^2 +
\lambda\|w_\e\|_{L^2}^2\\
&&\qquad\le \limsup_{\e\to0}
\|w_\e\|_{\calD^{s,2}}^2 +
\lambda\|w_\e\|_{L^2}^2\le \limsup_{\e\to0}\nu(V_\e)=\ti{\nu}.\end{eqnarray*}
This gives that
$$ \lim_{\e\to0}
\|w_\e\|_{\calD^{s,2}}^2 +
\lambda\|w_\e\|_{L^2}^2 =
\|w\|_{\calD^{s,2}}^2 +
\lambda\|w\|_{L^2}^2.$$
By making use of this and of the weak convergence of~$w_\e$, we infer that~$w_\e\to w$
in the Hilbert norm $\sqrt{ \|\cdot\|_{\calD^{s,2}}^2 +
\lambda\|\cdot\|_{L^2}^2}$.

Since this norm is equivalent to the one in~$H^s$, we have proved
that~$w_\e\to w $ in $H^s$.~\QED

\section{Completion of the proof of Theorem~\ref{TH:CON:2}}\label{CCC0}

Now we finish the proof of Theorem~\ref{TH:CON:2}. For this,
we suppose that~$V$ has a unique global minimum point at~$x_0$
Let $u_\e$ be a minimizer for $\nu_\e(V)$. Then $v_\e (x):= u_\e(x_0+\e x)$ is a minimizer for $\nu(V_\e)$.
By Lemma \ref{no0l}, there are     points $a_\e, c\in\R^N$ such that, up to a subsequence,
$$
w_\e(x):=v_\e(x+a_\e)=\ti{U}(x- c)+\o_{\e}(x),
$$
where  $\|\o_\e\|_{H^{s}}\to 0$, the function~$\ti{U}$ is a minimizer
for~$\nu(\ti{V})$, and, comparing~\eqref{DEFTIV}
and~\eqref{DEFTIV-1}, we have that
\begin{equation}\label{Cop}
\lim_{\e\to 0} V(\e a_\e+x_0)=\ti{V}=\min_{x\in \R^N} V(x)=V(x_0).\end{equation}
Now we prove that
\begin{equation}\label{Cop2}
\lim_{\e\to0}\e a_\e=0.\end{equation}
Suppose not, say~$|\e a_\e|\ge a_0$ for some~$a_0>0$ and an infinitesimal
sequence of~$\e$'s. Then $|\e a_\e|$ remains bounded,
otherwise, by~\eqref{GHJ}, the limit in~\eqref{Cop}
would be strictly larger than~$V(x_0)$.

Accordingly, there exists
an infinitesimal
sequence of~$\e$'s for which~$\e a_\e\to \alpha$, for some~$\alpha\in\R^N$
with~$|\alpha|\ge a_0>0$. {F}rom this and~\eqref{Cop},
we obtain that
$$ V(x_0)=\lim_{\e\to 0} V(\e a_\e+x_0)=V(\alpha+x_0).$$
This contradicts the
uniqueness of the minimal point for~$V$,
and so it proves~\eqref{Cop2}.

Now we claim that
\be\label{eq:Copactinfty}
\sup_{\e}\int_{|x|\geq R } w_\e^{r }\, dx\to 0  \quad \textrm{ as } R\to  \infty,
\ee
with\footnote{If~$N\le 2s$, the above definition of~$r$
can be replaced by just fixing~$r\in (1,+\infty)$.}
$r:=\frac{2N}{N-2s}$.
To see this, we can assume by contradiction that there
exists $\d$ positive and a sequence of $R_n\to \infty$ such that
$$
\sup_{\e}\int_{|x|\geq R_n } w_\e^{r }\, dx \geq \d  \quad \textrm{ as } n\to  \infty,
$$ 
This implies that for a sequence of  $\e_n \to 0$, we have 
$$
 \int_{|x|\geq R_n } w_{\e_n}^{r }\, dx \geq \d  \quad \textrm{ as } n\to  \infty.
$$
Because $w_{\e_n}$ converges strongly in $L^r$, we have
(see e.g.~\cite[Theorem 4.9]{Brezis})
that there exists $h\in L^r$ and a  subsequence, still denoted by $\e_n$ such that $w_{\e_n}\leq h$ a.e. in $\R^N$. But then 
$$
0< \d \leq  \int_{|x|\geq R_n } w_{\e_n}^{r }\, dx\leq   \int_{|x|\geq R_n } h^{r }\, dx \to 0 \textrm{ as } n\to  \infty.
$$
This leads to a contradiction. We thus have proved \eqref{eq:Copactinfty}.

Next we observe that $\Ds w_\e -\nu(\ti{V}_\e) w_\e^{p-1} w_\e\leq 0$ in $\R^N$. Since $w_\e^{p-1}\in L^q_{loc}$ for some 
$q>\frac{N}{2s}$, we deduce from \cite[Proposition 2.6]{JLX}
  that for any compact set $K$, we have
$$
\max_{K} w_\e \leq C  \int_{K } w_\e^{r }\, dx,
$$
where $r$ is as above. 
We therefore conclude from \eqref{eq:Copactinfty} that 
$$
\sup_{\e} w_\e(x)\to 0 \quad  \textrm{ as } |x|\to \infty.
$$
This together with Lemma C.2 in~\cite{FLS} also imply that 
\be\label{eq:weleqmxnp2s}
w_\e(x)\leq \frac{Const}{1+|x|^{N+2s}}.
\ee
By scaling back, we obtain
\be\label{eq:ueconcentx0}
 u_\e(x)=v_\e\left( \frac{x-x_0}{\e}\right)=w_\e\left( \frac{x-x_0-\e a_\e}{\e}\right)
\le \frac{Const\;\e^{N+2s}}{\e^{N+2s}+|x-x_0-\e a_\e|^{N+2s}}.
\ee
It is then clear that $u_\e$ concentrates at $x_0$ in the sense of~\eqref{eq:concent:weak}.\\

 Now to prove the last statement of the theorem ($u_\e$ has a unique global maximum point), we observe that $u_\e\in C^{2,\a}_{loc}$ and by \eqref{eq:ueconcentx0}, we have  $\lim\limits_{|x|\to\infty}u_\e(x)=0$ for every fixed and positive $\e$. We can therefore let  $u(x_\e)=\max\limits_{\R^N}u_\e$.   Then $\Ds u_\e(x_\e)\geq 0$ and thus from \eqref{MAIN:EQ} (recalling \eqref{eq:Ass-V}),   we deduce that 
$$
u(x_\e)\geq \left(\frac{\bar V}{\nu(V)} \right)^{\frac{1}{p-1}}=: C_0.
$$
Hence by \eqref{eq:ueconcentx0}, we get 
$$
C_0 \leq  \frac{Const\;\e^{N+2s}}{\e^{N+2s}+|x_\e-x_0-\e a_\e|^{N+2s}}
$$
so that 
\be \label{eq:xe-x0eaeleqC1e}
|x_\e-x_0-\e a_\e|\leq C_1 \e.
\ee 
From this we conclude, provided $|x-x_\e|\geq \e R\geq 2\e C_1$, that 
$$
u_\e(x)\leq \frac{Const}{1+R-{C_1}} \leq \frac{Const}{1+R/2}
$$
and this completes the proof of concentration of $u_\e$ at $x_0$.\\
We now prove the uniqueness of $x_\e$.  Indeed,  we observe that 
\begin{align*}
\Ds \o_\e=\Ds w_\e-\Ds\ti{U}(\cdot-c)&=V(\e x+\e a_\e+x_0) [\ti{U}(\cdot-c) -w_\e]+[V(x_0)- V(\e x+\e a_\e+x_0)]\ti{U}(\cdot-c)\\
&+
[\nu(V_\e)-\nu(V(x_0))] w_\e^p+ \nu(V(x_0))[ w_\e^p  -\ti{U}^p(\cdot-c)] .
\end{align*}
We rewrite this as
$$
\Ds \o_\e+\b_\e(x)\o_\e= [V(x_0)- V(\e x+\e a_\e+x_0)]\ti{U}(\cdot-c)+[\nu(V_\e)-\nu(V(x_0))] w_\e^p,
$$
where we have set  
$$ \b_\e(x) =V(\e x+\e a_\e+x_0)-\frac{\nu(V(x_0))[ w_\e^p  -\ti{U}^p(\cdot-c)]}{w_\e-\ti{U}(\cdot-c)}.$$
By \eqref{eq:weleqmxnp2s}, we have $ | w_\e^p  -\ti{U}^p(\cdot-c)|\leq C |w_\e-\ti{U}(\cdot-c)|$ and thus $ |\b_\e(x)|\leq Const.$. Applying \cite[Proposition 2.6]{JLX}, we deduce that $\o_\e \to 0$ in $C^{0,\a}_{loc}(\R^N)$ for some $\a\in (0,1)$. Now by a bootstrap argument and using  Proposition 2.1.8 in \cite{Sil}, we conclude that $\o_\e=w_\e-\ti{U}(\cdot-c)\to 0 $ in   $C^{2,\a}_{loc}(\R^N)$ for some $\a\in (0,1)$.\\
We now set $\bar{w}_\e(x)=w_\e(x+\bar{x}_\e)$ with $\bar{x}_\e=\frac{x_\e-x_0-\e a_\e}{\e}$. We notice that 0 is the global maximum point 
of $\bar{w}_\e$ and so we have 
$$
0=\n\bar{w}_\e(0)=\n \ti{U}(\bar{x}_\e-c)+\n \o_\e(\bar{x}_\e).
$$
Recalling that $\ti{U}$ is symmetric decreasing with respect to
the origin, that has
a unique critical point and also $\o_\e\to 0\in C^{2,\a}_{loc}(\R^N)$. Therefore from  \eqref{eq:xe-x0eaeleqC1e} we deduce that 
$$
|\bar{x}_\e-c|\to 0.
$$
It is clear that any other global maximum point of $\bar{w}_\e$ must stay in a neighborhood of $c$. We then observe that 
$$
\bar{w}_\e(x)=\ti{U}(x)+\bar{\o}_\e(x),
$$
where $ \bar{\o}_\e(x)=[\ti{U}(x+\bar{x}_\e-c)-\ti{U}(x) ]+{\o}_\e(x)$. Since $\ti{U}\in C^\infty(\R^N)$, 
we obtain $\bar{\o}_\e\to0\in C^{2,\a}_{loc}$.  Now using  Lemma 4.2 in \cite{NiTa}, we conclude that the only critical point for $\bar{w}_\e$ is the origin.
\QED
\begin{Remark}\label{rem:uniquxe}
{\rm We remark that from the above proof, the minimizers $u_\e$ for $\nu_\e(V)$ has the following precise form:
$$
u_\e(\e x+x_\e)=\ti{U}(x)+\bar{\o}_\e(x),
$$
where $\bar{\o}_\e\to0$ in $H^s(\R^N)\cap C^{2,\a}_{loc}(\R^N)\cap L^\infty(\R^N)$ with $x_\e$ the unique global maximizer for $u_\e$ and $x_\e$ converges to $x_0$
which is the global minimum point for $V$. Also $\ti{U}$ is the unique minimizer for $\nu(V(x_0))$.}
\end{Remark}

\section{Non-degeneracy and uniqueness: preliminaries for
the proof of Theorem~\ref{th:uniq-rad}}\label{PR1}

Now we will deal with the functional
\begin{equation}\label{P17} J_\e(u,\nu(V_\e)):=\frac{1}{2}\|u\|_\e^2-\frac{\nu(V_\e)}{p+1}
\int_{\R^N}|u|^{p+1}\,dx\end{equation}
and we will consider
the scalar products that induce the
norms of the fractional spaces used in this paper, namely we set
\begin{eqnarray*}
&& \la u,v\ra_{\calD^{s,2}}:=
\int_{\R^N} |\xi|^{2s} \widehat{u}\,\widehat{v}\,d\xi,\\
&& \la u,v\ra_{\e,V}:=\int_{\R^N} V(\e x+x_0) \, u(x)\,v(x)\,dx,\\
{\mbox{and }} &&
\la u ,v\ra_\e := \la u,v\ra_{\calD^{s,2}} + \la u,v\ra_{\e,V}.
\end{eqnarray*}
The Hilbert space associated with~$\la \cdot,\cdot\ra_\e$
will be denoted by~$H^s_\e$ and, as usual, we say that~$u
\perp_\e v$ whenever~$\la u ,v\ra_\e=0$.
One simple, but important feature, is that the
radially symmetric minimizer~$U$
for $\nu(V_0)$ is perpendicular in~$H^s_0$ (that is~$H^s_\e$
with~$\e=0$) to its derivatives,
and the derivatives themselves are perpendicular to each other,
according to the following result:

\begin{Lemma}\label{ij}
For any~$i\in\{1,\dots,N\}$, we have that
\be\label{first}\begin{split}
&\la U ,\partial_i U\ra_0=0\\
{\mbox{and }} \ & \int_{\R^N} U^{p}\partial_i U=0.\end{split}\ee
Moreover, for any~$i$, $j\in\{1,\dots,N\}$, with~$i\ne j$,
we have that
\be\label{es0}
\int_{\R^N} U^{p-1}\,\partial_i U\,\partial_j U
=0\ee
and
\be\label{sec}\la \partial_i
U ,\partial_j U\ra_0=0.\ee
\end{Lemma}

\proof By construction
\be\label{EU0}
(-\Delta)^s U+V(0)U=U^p
\ee
and so, taking derivatives,
\be\label{EUi}
(-\Delta)^s (\partial_i U)+V(0)\partial_i U=p U^{p-1} \partial_i U.
\ee
We multiply~\eqref{EU0} by~$\partial_i U$ and~\eqref{EUi} by~$U$ and
integrate: we obtain, respectively,
$$ \la U ,\partial_i U\ra_0 = \int_{\R^N} U^p\,\partial_i U
\ {\mbox{ and }} \
\la U ,\partial_i U\ra_0 = p\int_{\R^N} U^p\,\partial_i U
.$$
By comparing these two equations we obtain that
$$ p\int_{\R^N} U^p\,\partial_i U=\int_{\R^N} U^p\,\partial_i U,$$
and so
$$ \la U ,\partial_i U\ra_0 = \int_{\R^N} U^p\,\partial_i U=0,$$
that proves~\eqref{first}.

Now we use the rotational invariance of~$U$
to write~$U(x)=\bar U(|x|)$, for some~$\bar U:\R\to\R$.
Then we have that~$\partial_i U(x) =\bar U'(|x|) \,|x|^{-1}\,x_i$
and so, by symmetry
$$ \int_{\R^N} U^{p-1}(x)\,\partial_i U(x)\,\partial_j U(x)\,dx
=\int_{\R^N} \bar U^{p-1}(|x|)\;\big|\bar U'(|x|)\big|^2\,|x|^{-2} \,x_i x_j
\,dx =0.$$
This establishes~\eqref{es0}.
Then, formula~\eqref{sec} follows multiplying~\eqref{EUi}
by~$\partial_j U$ and integrating over~$\R^N$.~\QED
Our next result is of coercivity type. It  is stronger than what we will need in the following of the paper we expose it here because we believe that it might be of interest.\\

We also mention that in the rest of the paper, the regularity assumption   can be relaxed to $V\in C^1(\R^N)$.\\

Given the radially symmetric minimizer~$U$ for $\nu(V_0)$ and~$a
\in \R^N$, we define $U_a(x):=U(x-a)$ and
\be\label{8csdfww87}
W_\e:=\Big\{ v\in H^{s} {\mbox{ s.t. $v\perp_\e U_a$ and
$v\perp_\e \de_j U_a\,$ for any $\,j=1,\dots,N$}}\Big\}.\ee
With this, we can bound the second derivative of~$J_\e(U_a,\nu_\e)$
from below as follows:

\begin{Lemma}\label{lem:coerciv}
Let~$J_\e$ be as in~\eqref{P17}. There exists~$\e_0>0$ such that for any~$\e
\in(0,\e_0)$, for any~$v\in W_\e$ and for any~$a\in\R^N$
$$
J_\e''(U_a,\nu(V_\e))[v,v]\geq \,Const\, \|v\|^2_{\e}. $$
The $\,Const\,$ above does not depend on~$a$.
\end{Lemma}

\proof Up to translations, we can suppose that~$x_0=0$.
We consider $\chi\in C^\infty_c(\R^N,(0,2))$ such that $\chi=1$ in
$B_1$ and~$\chi=0$ in $\R^N\setminus B_2$.
Also we take~$R>1$, to be chosen suitably large in the sequel.
We define
\begin{eqnarray*}
&& \chi_R(x)=\chi\left(a+\frac{x}{R}\right),\\
&& \bar\chi_R:=1-\chi_R,\\
&& v_1:=\chi_R v,\\
&& v_2:=\bar\chi_R v,\\
{\mbox{and }} && I_1:= \int_{\R^{2N}} \frac{\chi_R(x)\,\bar\chi_R(x)\,
(v(x)-v(y))^2}{|x-y|^{N+2s}}\,dx\,dy.
\end{eqnarray*}
First we prove that
\begin{equation}\label{-01-}
\big| \la v_1,v_2\ra_{\calD^{s,2}}\big|\le \eta_R(v),
\end{equation}
with~$\eta_R(v)$ not depending on~$\e$ and such that
\be\label{ETA} \lim_{R\to+\infty}\eta_R(v)=0,\ee
for $v$ fixed.
To this goal, we compute
\begin{eqnarray*}
&&(v_1(x)-v_1(y))\,(v_2(x)-v_2(y))\\
&=& \Big( v(x)\chi_R(x)-v(y)\chi_R(y)\Big)\,
\Big( v(x)\bar\chi_R(x)-v(y)\bar\chi_R(y)\Big) \\
&=&\left( v(x)\Big(\chi_R(x)-\chi_R(y)\Big)
+\chi_R(y)\Big(v(x)-v(y)\Big)\right)\\
&& \quad\cdot\,
\left( \bar\chi_R(x)\Big(v(x)-v(y)\Big)+
v(y)\Big( \bar\chi_R(x)-\bar\chi_R(y)\Big)\right)
\\ &=&
\left( v(x)\Big(\chi_R(x)-\chi_R(y)\Big)
+\chi_R(y)\Big(v(x)-v(y)\Big)\right)\\&& \quad\cdot\,
\left( \bar\chi_R(x)\Big(v(x)-v(y)\Big)-
v(y)\Big( \chi_R(x)-\chi_R(y)\Big)\right) \\
&=&- v(x)\,v(y)\Big(\chi_R(x)-\chi_R(y)\Big)^2+
\chi_R(x)\,\bar\chi_R(x) \Big(v(x)-v(y)\Big)^2\\
&&\quad
+v(x)\,\bar\chi_R(x) \Big(v(x)-v(y)\Big)\Big(\chi_R(x)-\chi_R(y)\Big)\\
&&\quad
-v(y)\,\chi_R(y) \Big(v(x)-v(y)\Big)\Big(\chi_R(x)-\chi_R(y)\Big).
\end{eqnarray*}
Therefore
\be\label{4.1}
\big|\la v_1,v_2\ra_{\calD^{s,2}}\big|\,\le\, I_1\,+\,Const\,(J_1+J_2)
\ee
with
\begin{equation}\label{dmu}\begin{split}
J_1 \,:=& \int_{\R^{2N}} |v(x)|\,|v(y)|\,
\Big(\chi_R(x)-\chi_R(y)\Big)^2
\,d\mu(x,y),\\
J_2 \,:=& \int_{\R^{2N}} |v(x)|\;|v(x)-v(y)|\;
|\chi_R(x)-\chi_R(y)|
\,d\mu(x,y),\\
{\mbox{and }} \
d\mu(x,y) \,:=& |x-y|^{-N-2s}\,dx\,dy.
\end{split}\end{equation}
Now we observe that~$\|\nabla\chi_R\|_{L^\infty}\le \,Const\, R^{-1}$, and so
\begin{equation}\label{sym0}
|\chi_R(x)-\chi_R(y)|\le \,Const\,\min\big\{ 1, \ R^{-1} |x-y|\big\}.
\end{equation}
Therefore, for any~$x\in\R^N$,
\begin{equation*}\begin{split}
&\int_{\R^n} \frac{|\chi_R(x)-\chi_R(y)|^2}{|x-y|^{N+2s}}\,dy \\ \le\,&
\,Const\,\left[
\int_{B(0,R)} \frac{R^{-2}\,|x-y|^2}{|x-y|^{N+2s}}\,dy +
\int_{\R^N\setminus B(0,R)} \frac{1}{|x-y|^{N+2s}}\,dy
\right] \\
=\,& \,Const\,R^{-2s}.
\end{split}\end{equation*}
Using this and
the H\"older inequality we obtain
\begin{equation}\label{4.1.1}\begin{split}
J_2 \,&\le \sqrt{\int_{\R^{2N}} |v(x)|^2\,|\chi_R(x)-\chi_R(y)|^2\,d\mu(x,y)}\,\cdot\,\sqrt{\int_{\R^{2N}} |v(x)-v(y)|^2\,d\mu(x,y)}
\\ &\le \sqrt{Const\,R^{-2s}\int_{\R^N}|v(x)|^2\,dx}\,\cdot\,\|v\|_{\calD^{s,2}}
\\ &\le\, Const\,R^{-s} \,\|v\|_{L^2}\,\cdot\,\|v\|_{\calD^{s,2}}
\\ &\le\, Const\,R^{-s} \,\|v\|_\e^2.
\end{split}\end{equation}
Now we define
\begin{eqnarray*}
&& \R^{2N}_R := \big\{ (x,y)\in\R^{2N} {\mbox{ s.t. }} |x-y|<R\big\}\\
{\mbox{and }}&& {\mathcal{V}}:=
\big\{ (x,y)\in\R^{2N} {\mbox{ s.t. }} |v(x)|\ge |v(y)|\big\}.
\end{eqnarray*}
By symmetry
\begin{equation}\label{sym1}
\begin{split}
& \int_{\R^{2N}_R} |v(x)|\,|v(y)|\,|x-y|^2\,d\mu(x,y)\le
2\int_{\R^{2N}_R\cap {\mathcal{V}} } |v(x)|^2\,|x-y|^2\,d\mu(x,y) \\
&\qquad \le 2\int_{\R^N} |v(x)|^2 \,\left[ \int_{B(x,R)} |x-y|^{2-N-2s}\,dy
\right]\,dx
\le \,Const\,R^{2-2s} \int_{\R^N} |v(x)|^2\,dx\\
&\qquad= \,Const\,R^{2-2s}\,\|v\|_{L^2}^2.
\end{split}
\end{equation}
Similarly,
\begin{equation*}
\begin{split}
& \int_{(\R^{2N}\setminus\R^{2N}_R)} |v(x)|\,|v(y)|\,d\mu(x,y)\le
2\int_{(\R^{2N}\setminus\R^{2N}_R)\cap {\mathcal{V}} } |v(x)|^2\,d\mu(x,y) \\
&\qquad \le 2\int_{\R^N} |v(x)|^2 \,\left[ \int_{\R^N\setminus B(0,R)} |x-y|^{-N-2s}\,dy\right]\,dx
\le \,Const\,R^{-2s} \int_{\R^N} |v(x)|^2\,dx\\
&\qquad= \,Const\,R^{-2s}\,\|v\|_{L^2}^2.
\end{split}
\end{equation*}
We use the latter inequality together with~\eqref{sym0} and \eqref{sym1}
to conclude that
\begin{eqnarray*}
J_1 &\le& Const\,\left[R^{-2}\int_{\R^{2N}_R} |v(x)|\,|v(y)|\,|x-y|^2\,d\mu(x,y)+
\int_{\R^{2N}\setminus\R^{2N}_R} |v(x)|\,|v(y)|\,d\mu(x,y)\right]\\
&\le& Const\,R^{-2s}\,\|v\|_{L^2}^2.
\end{eqnarray*}
Hence, by~\eqref{4.1} and~\eqref{4.1.1},
\be\label{9i08} \big|\la v_1,v_2\ra_{\calD^{s,2}}\big|\le I_1+\,Const\,(R^{-2s}+R^{-s})\|v\|_\e^2.\ee
Now we estimate~$I_1$.
For this we observe that the function~$\chi_R\bar\chi_R$ is supported
in~$B(a,2R)\setminus B(a,R)$, hence
$$ I_1\le \int_{(B(a,2R)\setminus B(a,R))\times \R^{N}} \frac{
(v(x)-v(y))^2}{|x-y|^{N+2s}}\,dx\,dy.$$
Since~$v$ is a fixed function of~$H^s$, we have that
$$ \lim_{R\to+\infty} \int_{(B(a,2R)\setminus B(a,R))\times \R^{N}} \frac{
(v(x)-v(y))^2}{|x-y|^{N+2s}}\,dx\,dy=0.$$
These considerations and~\eqref{9i08} imply~\eqref{-01-}, as desired.

{F}rom~\eqref{-01-}, we obtain that
\begin{equation}\label{ext D}\begin{split}
& \|v\|_{\calD^{s,2}}^2=\|v_1+v_2\|_{\calD^{s,2}}^2=
\|v_1\|_{\calD^{s,2}}^2+\|v_2\|_{\calD^{s,2}}^2+2\la v_1,v_2\ra_{\calD^{s,2}}
\\ &\qquad\le \|v_1\|_{\calD^{s,2}}+\|v_2\|_{\calD^{s,2}}
+2\eta_R(v).\end{split}
\end{equation}
Moreover
\begin{equation*}
\begin{split}
 \|v\|^2_{\e,V}\,&
= \|v_1\|^2_{\e,V}+
 \|v_2\|^2_{\e,V}+\int_{\R^N} V(\e x) \, v_1(x)\,v_2(x)\,dx
\\&\le \,Const\,\big(  \|v_1\|^2_{\e,V}+
 \|v_2\|^2_{\e,V}\big).
\end{split}\end{equation*}
This and~\eqref{ext D} yield that
\begin{equation}\label{J7E}
 \|v\|^2_\e\le \,Const\,\big(\|v_1\|^2_\e+
 \|v_2\|^2_\e +\eta_R(v)\big).
\end{equation}
On the other hand, $v_1 v_2=\chi_R \bar\chi_R v^2$, therefore
$v_1v_2\ge0$ and it is supported in~$B(a,2R)\setminus B(a,R)$.
In this domain~$U_a$ is of the order ~$R^{-(N+2s)}$, therefore
$$ \int_{\R^N} U_a^{p-1} v_1\,v_2\le \,Const\, R^{-(p-1)(N+2s)}
\int_{B(a,2R)\setminus B(a,R)}|v|^2
\le \,Const\, R^{-(p-1)(N+2s)}\,\|v\|^2_{L^2}.$$
{F}rom this and \eqref{-01-} we infer that
\begin{equation}\label{eq:Jv1v2}\begin{split}
J_\e''(U_a,\nu(V_\e))[v_1,v_2]&= \la v_1,v_2\ra_{\calD^{s,2}}+
\int_{\R^N}V(\e x)v_1v_2 -p\nu(V_\e)\int_{\R^N}U_a^{p-1}v_1v_2\\
&\ge \,-\,Const\,R^{-s}\,\|v\|_\e^2+0-\,Const\, R^{-(p-1)(N+2s)}\,\|v\|^2_{L^2}\\
&\ge \,-\,Const\,R^{-\gamma}\,\|v\|_\e^2,
\end{split}\end{equation}
up to renaming constants, where~$\gamma:=\min\{s,\,(p-1)(N+2s)\}>0$ (here we have also used Lemma~\ref{no0l}
to bound~$\nu(V_\e)$ uniformly in~$\e$). Similarly, $v_2$ is supported
outside~$B(0,R)$, hence
$$ \int_{\R^N}U_a^{p-1}v_2^2 \le\,Const\,R^{-(p-1)(N+2s)} \int_{\R^N} v^2,$$
and therefore
\begin{align}\label{eq:Jv2v2}
J_\e''(U_a,\nu(V_\e))[v_2,v_2]&= \|v_2\|^2_\e-p \nu(V_\e)
\int_{\R^N}U_a^{p-1}v_2^2\geq \|v_2\|^2_\e - \,Const\,
R^{-(p-1)(N+2s)}\,\|v\|^2_{L^2}.
\end{align}
Next we estimate $ J_\e''(U_a,\nu(V_\e))[v_1,v_1]$. To this goal, we
project~$v_1$ along the space spanned by~$U_a$
and its derivatives, i.e. we set
$$
\psi:=\frac{1}{\|U_a\|^2_0}\la v_1,U_a\ra_0 U_a+
\frac{1}{\|\de_iU_a\|^2_0}\la v_1,\de_iU_a\ra_0 \de_i U_a,
$$
where the repeated indices convention is used,
and~$w:=v_1-\psi$.
Therefore
\be
J_\e''(U_a,\nu(V_\e))[v_1,v_1]=J_\e''(U_a
,\nu(V_\e))[w,w]+J_\e''(U_a,\nu(V_\e))[\psi,\psi]+2 J_\e''(U_a,\nu(V_\e))[w,\psi].
\ee
We observe that the norms~$\|\cdot\|_0$
and~$\|\cdot\|_\e$ are comparable, thanks to~\eqref{eq:Ass-V}. Therefore
\begin{equation}\label{ep}
|\psi| \,\le \frac{\|v_1\|_0}{\|U\|_0}\,U_a+
\frac{\|v_1\|_0}{\|\de_iU_a\|_0}\, |\de_i U_a| \le Const\, \|v_1\|_\e\,\Big( U_a+|\de_i U_a|\Big),
\end{equation}
hence
$$ \int_{\R^N} (1+|x|)\psi^2 \le \,Const\, \|v_1\|_\e^2.$$
Using this, the fact that~$|V(\e x)-V(0)|\le\,Const\,\e\,|x|$,
and that~$v_1$ is supported in~$B(a,R)$,
we conclude that
\begin{equation}\label{poe}
\begin{split}
&\int_{\R^N}|V(\e x)-V(0)|w^2\\
&= \int_{B(a,R)}|V(\e x)-V(0)|\,v_1^2+
\int_{\R^N}|V(\e x)-V(0)|\,\psi^2-2\int_{\R^N}|V(\e x)-V(0)|\,v_1\psi\\
&\le \,Const\,\left[\int_{B(a,R)}|V(\e x)-V(0)|\,v_1^2+
\int_{\R^N}|V(\e x)-V(0)|\,\psi^2\right]\\
&\leq \,Const\,\e \,(R+|a|)\, \|v_1\|_\e^2.
\end{split}
\end{equation}
Now we remark that~$w$ is orthogonal in~$H^s_0$ (i.e. in~$H^s_\e$
with~$\e=0$) to any element of the basis~$\{U_a,\,\partial_1 U_a,\dots,\,\partial_N U_a\}$, thanks to Lemma~\ref{ij}. Hence, from~\cite{FLS},
we have that
$$ J_0''(U_a,\nu_0)[w,w]\ge \,Const\,\|w\|^2_0\ge
\,Const\,\|w\|^2_\e.$$
As a consequence,
\be\label{eq:Jww}
\begin{split}
J_\e''(U_a,\nu(V_\e))[w,w]&= J_0''(U_a,\nu(V_0))[w,w]+
\int_{\R^N}[V(\e x)-V(0)]w^2-p( \nu(V_\e)-\nu(V_0))
\int_{\R^N}U_a^{p-1}w^2\\
&\geq \,Const\, \|w\|^2_\e+  \int_{\R^N}[V(\e x)-V(0)]w^2-p( \nu(V_\e)-\nu(V_0))
\int_{\R^N}U^{p-1}w^2\\
&\geq \,Const\, \|w\|^2_\e
-\,Const\,\e \,(R+|a|)\, \|v_1\|_\e^2-\,Const\,|\nu(V_\e)-\nu(V_0)|\,\|v_1\|_\e,
\end{split}\ee
where both~\eqref{poe}
and~\eqref{ep} were used in the last inequality.

Furthermore, since~$v\perp_\e U_a$, we have
\begin{align*}
\la v_1,U_a\ra_0&=\la v_1,U_a\ra_\e+\int_{\R^N}[V(0)-V(\e x)]v_1
U_a\\
&=\la v,U_a\ra_\e- \la v_2,U_a\ra_\e  +\int_{\R^N}[V(0)-V(\e x)]v_1
U_a\\
&=-  \la v_2,U_a\ra_\e  +\int_{\R^N}[V(0)-V(\e x)]v_1 U_a\\
&= - \la v_2,U_a\ra_{\calD^{s,2}}- \int_{\R^N}V(\e x)v_2 U_a
+\int_{\R^N}[V(0)-V(\e x)]v_1 U_a\\
&= -\int_{\R^N}  v_2\Ds U_a- \int_{\R^N}V(\e x)v_2 U_a
+\int_{\R^N}[V(0)-V(\e x)]v_1 U_a\\
&=-\int_{\R^N}  v_2[-V(0) U_a+ \nu_0 U_a^p]- \int_{\R^N}V(\e x)v_2
U_a+\int_{\R^N}[V(0)-V(\e x)]v_1 U_a \\
&= \int_{\R^N}[V(0)-V(\e x)]v U_a- \nu_0\int_{\R^N} U_a^p v_2\end{align*}
thus, since~$v_2$ is supported outside~$B(a,R)$,
\be\label{q1q}
\Big| \la v_1,U_a\ra_0\Big|
\le \,Const\,\Big(\e\,(R+|a|)\,\|v\|_{L^2}+ R^{-(p-1)(N+2s)}\|v\|_{L^2}\Big).\ee
In a similar way, since also~$v\perp_\e\partial_i U_a$,
we have that
\be\label{q2q}\Big|\la v_1,\de_iU_a\ra_0\Big|
\le \,Const\,\Big(\e\,(R+|a|)\,\|v\|_{L^2}+ R^{-(p-1)(N+2s)}\|v\|_{L^2}\Big).\ee
We deduce from~\eqref{q1q} and~\eqref{q2q} that
$$\|\psi\|_0\le \,Const\,\Big( \Big|\la v_1,U_a\ra_0\Big|+
\Big|\la v_1,\de_iU_a\ra_0\Big|\Big)\le
\,Const\,\Big(\e(R+|a|)\,\|v\|_{L^2}+ R^{-(p-1)(N+2s)}\|v\|_{L^2}\Big)$$
and so, since the two norms are comparable,
$$\|\psi\|_\e
\le \,Const\,\Big(\e(R+|a|)\,\|v\|_{L^2}+ R^{-(p-1)(N+2s)}\|v\|_{L^2}\Big).$$
So, we use the fact that
$$ 2|\la v_1,\psi\ra_\e|=2|\la v_1/2,\, 2\psi\ra_\e|\le
\frac{\|v_1\|_\e}{4}+4\|\psi\|_\e$$
to conclude that
\begin{eqnarray*}
\|w\|_\e^2 &=& \|v_1\|^2_\e+\|\psi\|^2_\e-2\la v_1,\psi\ra_\e
\\ &\ge& \frac{3}{4} \|v_1\|^2_\e - \,Const\,
(\e(R+|a|)+ R^{-(p-1)(N+2s)})^2\,\|v\|_{L^2}^2.
\end{eqnarray*}
Exploiting this and~\eqref{eq:Jww} we obtain
\be\label{eq:Jv1v1} \begin{split} & J_\e''(U_a,\nu(V_\e))[w,w]\ge
\,Const\, \|v_1\|^2_\e - \,Const\,
(\e(R+|a|)+ R^{-(p-1)(N+2s)})^2\,\|v\|_{L^2}^2\\ &\qquad\quad
-\,Const\,\e \,(R+|a|)\, \|v_1\|_\e^2-\,Const\,|\nu(V_\e)-\nu(V_0)
|\,\|v_1\|_\e.\end{split}\ee
Notice now that
$$
J_\e''(U_a,\nu(V_\e))[v,v]=J_\e''(U_a,\nu(V_\e)
)[v_1,v_1]+J_\e''(U_a,\nu(V_\e))[v_2,v_2]+2 J_\e''(U_a,\nu(V_\e))[v_1,v_2].
$$
Thus, by collecting \eqref{eq:Jv1v2}, \eqref{eq:Jv2v2} and \eqref{eq:Jv1v1},
we obtain
\begin{align*}
J_\e''(U_a,\nu(V_\e))[v,v]&\geq  \, Const\,
(\|v_1\|^2_\e+\|v_2\|^2_\e )
\\ &\qquad\,-\,Const\, (\e(R+|a|)+R^{-\gamma} )\,\|v\|_\e^2,
-\,Const\,
R^{-(p-1)(N+2s)}\,\|v\|^2_{L^2}
\\&\qquad-\,Const\,\e \,(R+|a|)\, \|v_1\|_\e^2-
\,Const\,|\nu(V_\e)-\nu(V_0)|\,\|v_1\|_\e
\end{align*} Now, recalling~\eqref{J7E}
and~\eqref{ETA},
and sending first~$\e\to 0$ and then~$R\to +\infty$, we get the
desired result.~\QED

\section{Uniqueness of radial solutions}
In this section we assume that $V$ is  radial and we consider the functional
in~\eqref{P17}. We 
denote by $H^s_r$ the subspace of $H^s$ of radially symmetric function.
We will make use of the minimizer~$U$ for~$\nu(V_0)$,
normalized with~$\|U\|_{L^{p+1}}=1$, which is a solution of
\be\label{U weak}
\langle U,v\rangle_{\calD^{s,2}}
+V(0)\,\langle U,v\rangle_{L^2}=\nu(V_0)\,\int_{\R^N} U^p(x)\,v(x)\,dx,
\ee
for every~$v\in H^s$.

We also define
$I_\e$ as the restriction of
$u\mapsto  J_\e(u,\nu(V_\e))$ on $H^s_r$.
Next, we define the operator~$\Phi_{\e}: H^s_r\to H^{s}_r$
by \be\label{eq:defPhias} \Phi_{\e}(\o):=  I_{\e}'\left( U+\o\right). \ee

By \eqref{eq:defPhias}, we mean: for all $w\in H^s_r$ \be \la
\Phi_{\e}(\o),w\ra=  I_{\e}'\left( U+\o \right)[ w]. \ee

\begin{Lemma}\label{lem:invPHI}
There exist~$\delta,\e_0>0$ sufficiently small such that: for every $\e\in(0,\e)$,
if~$\Phi_\e(w_1)=\Phi_\e(w_2)$
for some~$w_1$, $w_2\in H^s_r$ with~$\|w_1\|_{\e}+\|w_2\|_{\e}\le\delta$,
then~$w_1=w_2$.
\end{Lemma}
\proof
The proof is a
consequence of  Lemma \ref{lem:coerciv}.
The details go as follows. First we fix the following notation: given~$f\in H^s_r$,
we define
$$ c_f:= \frac{\langle f,U\rangle_\e}{\| U \|_\e^2}
\ {\mbox{ and }} \ \ti{f}:= f-c_f U.$$
Notice that~$\ti{f}$ is radial, since so are~$f$ and~$U$,
and that~$\langle \ti{f},U\rangle_\e=0$. As a matter of fact,
since both~$\ti{f}$ and~$U$ are radial, a direct computation
based on odd symmetry shows that also~$\langle \ti{f},\partial_iU\rangle_\e=0$,
that is
\be\label{WE}
\ti{f}\in W_\e,\ee 
according to the definition in~\eqref{8csdfww87}.

Notice that~$f=\ti{f}+c_f U$.
We also consider the reflection of~$f$ with respect to~$U$, namely
\be\label{STAR}f^\star:=\ti{f}-c_f U.\ee
Now we observe that
\be\label{p-1 U}
J''_0(U,\nu(V_0)) [U,v]=
(1-p)\,\nu(V_0)\,\int_{\R^N} U^{p}(x)\,v(x)\,dx\ee
for any~$v\in H^s_r$. Indeed, for any~$v\in H^s$,
\begin{eqnarray*}
J''_0(U,\nu(V_0))[U,v]
&=&\langle U,v\rangle_{\calD^{s,2}}
+V(0)\,\langle U,v\rangle_{L^2}-p\nu(V_0)\,\int_{\R^N} U^{p-1}(x)\,U(x)\,v(x)\,dx
\\ &=&(1-p)\,\nu(V_0)\,\int_{\R^N} U^{p}(x)\,v(x)\,dx,
\end{eqnarray*}
thanks to~\eqref{U weak}, and this establishes~\eqref{p-1 U}.

Furthermore
$$ J''_\e(U,\nu(V_\e))[U,v]-J''_0(U,\nu(V_0))[U,v]=\int_{\R^N}
(V(\e x)-V(0))U(x)v(x)\,dx
-p(\nu(V_\e)-\nu(V_0))\int_{\R^N} U^{p}(x)\,v(x)\,dx.$$
This, combined with~\eqref{p-1 U} gives that
$$ J''_\e(U,\nu(V_\e))[U,v] =
\int_{\R^N}
(V(\e x)-V(0))U(x)v(x)\,dx
-c_\e\,\int_{\R^N} U^{p}(x)\,v(x)\,dx,$$
where
$$ c_\e:= -(1-p)\,\nu(V_0) +p(\nu(V_\e)-\nu(V_0)).$$
Notice that~$c_\e\to (p-1)\,\nu(V_0)>0$ as~$\e\to0$, due to
Lemma~\ref{no0l}.
In particular
$$ J''_\e(U,\nu(V_\e))[U,U] =\eta_\e
-c_\e,$$
with
$$ \eta_\e:=
\int_{\R^N}
(V(\e x)-V(0))U^2(x)\,dx \to 0,$$
as~$\e\to0$, by dominated convergence theorem.
We conclude that
\be\label{p-1 U-eps}
J''_\e(U,\nu(V_\e))[U,U] \le-\frac{c_\e}{2}\ee
for small~$\e$.
Now, for any $v$, $w\in H^s_r$, we set
\begin{eqnarray*}
\calN_{\e}(v)[w]&:=& 
\Phi_\e(v)[w]-\Phi_\e(0)[w]-\langle\Phi_\e'(0)[v],w\rangle\\
&=& I_{\e}'(U+v)[w]-I_{\e}'(U  )[w]-I_{\e}''(U )[v,w]\\
&=&\nu(V_\e)\left(-\int_{\R^N}|U+v|^pw\,
dx+\int_{\R^N}U^p w\,dx +p\int_{\R^N}U^{p-1}w\, dx
\right).
\end{eqnarray*}
Referring to page 128 in~\cite{AM}, we obtain
%%% $$|\calN_{\e}(v)[w]|\leq\,Const\,(\|v\|^2_\e+\|v\|^p_\e)\|w\|_\e $$ and
\be\label{78sdfghsee}
\|\calN_{\e}(v_1)-\calN_{\e}(v_2)\|\leq
\,Const\, (\|v_1\|_\e+
\|v_1\|^{p-1}_\e+\|v_2\|_\e+\|v_2\|^{p-1}_\e)\|v_1-v_2\|_\e
.\ee
Now we take~$w:=w_1-w_2$ and we use the notation in~\eqref{STAR}
and the assumption that~$\Phi_\e(w_1)=
\Phi_\e(w_2)$ to compute:
\begin{eqnarray*}
0 &=& \Phi_\e(w_1)[w^\star]-\Phi_\e(w_2)[w^\star]
\\ &=&
\calN_{\e}(w_1)[w^\star] +\Phi_\e(0)[w^\star]+
\langle\Phi_\e'(0)[w_1],w^\star\rangle
-
\calN_{\e}(w_2)[w^\star] -\Phi_\e(0)[w^\star]-
\langle\Phi_\e'(0)[w_2],w^\star\rangle\\
&=& \langle\Phi_\e'(0)[w_1],w^\star\rangle-
\langle\Phi_\e'(0)[w_2],w^\star\rangle
+\calN_{\e}(w_1)[w^\star]-\calN_{\e}(w_2)[w^\star]
\\ &=& J''_\e(U,\nu(V_\e)) [w,w^\star]
+\calN_{\e}(w_1)[w^\star]-\calN_{\e}(w_2)[w^\star].
\end{eqnarray*}
Thus, we write~$w=\ti{w}+c_w U$ and~$w^\star=\ti{w}-c_w U$, and we
exploit~\eqref{78sdfghsee} and~\eqref{p-1 U-eps}, to see that
\begin{eqnarray*}
0&\ge& J''_\e(U,\nu(V_\e)) [\ti{w},\ti{w}]-
c_w^2 J''_\e(U,\nu(V_\e)) [U,U]
-\,Const\, \delta^{\min\{1,p-1\}} \|w\|_\e \|w^\star\|_\e
\\ &\ge& J''_\e(U,\nu(V_\e)) [\ti{w},\ti{w}] +\frac{c_\e c_w^2}{2}
-\,Const\, \delta^{\min\{1,p-1\}} \|w\|_\e \|w^\star\|_\e.
\end{eqnarray*}
Now, thanks to~\eqref{WE}, we can make use of Lemma \ref{lem:coerciv}
and write that~$J_\e''(U,\nu(V_\e))[\ti{w},\ti{w}]\geq 
\,Const\, \|\ti{w}\|^2_{\e}$. So we obtain that
\be\label{0-11} 0\ge \,Const\, \Big(\|\ti{w}\|^2_{\e}+c_w^2\Big)-
\,Const\, \delta^{\min\{1,p-1\}} \|w\|_\e \|w^\star\|_\e.\ee
Also, by~\eqref{WE},
$$ \| w\|_\e^2 = \|\ti{w}\|_\e^2 + c_w^2\|U\|_\e^2+
2\langle \ti{w},U\rangle_\e=\|\ti{w}\|_\e^2 + c_w^2\|U\|_\e^2$$
and, similarly,
$$ \| w^\star\|_\e^2 
=\|\ti{w}\|_\e^2 + c_w^2\|U\|_\e^2.$$
In particular, $\| w\|_\e^2 \le\,Const\,\big( \|\ti{w}\|_\e^2 + c_w^2\big)$
and~\eqref{0-11} becomes
\begin{eqnarray*}
0&\ge& \,Const\, \Big(\|\ti{w}\|^2_{\e}+c_w^2\Big)-
\,Const\, \delta^{\min\{1,p-1\}} \|w\|_\e^2\\
&\ge& \big(\,Const\,-\,Const\,\delta^{\min\{1,p-1\}}\big)\|w\|_\e^2,
\end{eqnarray*}
which implies that~$\|w\|_\e=0$ if $\delta$ is small enough.
\QED

\section{Completeness of the proof of Theorem~\ref{th:uniq-rad}}\label{PR2}

Now we complete the proof of
Theorem~\ref{th:uniq-rad}.
For this, let $v_\e^i$ be a radial minimizer for $\nu(V_\e)$
in the class of radial competitors, with $i=1,2$. 
Then   by Lemma \ref{no0l}, provided $\e$ is sufficiently small, we have
$$
 v_\e^i(x)=  {U}+w_\e^i\qquad  \hbox{with}\quad  \|w^\e_i\|_\e\to0 \quad\hbox{as} \quad \e\to 0.
 $$ 
It turns out that $\Phi_\e(w_\e^i)=I_\e'(v_\e^i)=0$,
so we conclude that $w_\e^1=w_\e^2$,
due to Lemma~\ref{lem:invPHI}.
\QED

\section{Completeness of the proof of Corollary~\ref{coroll:uniq-rad}}\label{PR2bis}

Since $V$ is radial and radially
decreasing,
then using symmetric decreasing arguments, or the moving plane argument,
we have that $v_\e^i$ is radial.
Then the result follows by Theorem \ref{th:uniq-rad}.
\QED


\begin{thebibliography}{99}
 
\bibitem{AM} A. Ambrosetti and A. Malchiodi,
{\it Perturbation Methods and Semilinear Elliptic Problems on $\R^n$}.
Progress in Mathematics, Birkh\"auser Verlag, Basel-Boston-Berlin (2005).
 
\bibitem{Brezis} H. Brezis,
{\it Functional analysis, Sobolev spaces and partial differential equations}. 
Universitext. Springer, New York, 2011.

\bibitem{cabre-sire} X. Cabr\'{e} and Y. Sire,
{\it Nonlinear equations for fractional Laplacians, I:
Regularity, maximum principles, and Hamiltonian estimates}.
Ann. Inst. H. Poincar\'{e} Anal. Non Lin\'{e}aire 31 (2014), no. 1, 23-53.

\bibitem{CZ} G. Chen and Y. Zheng,
{\it Concentration phenomenon for fractional nonlinear Schr\"{o}dinger equations}.
Comm. Pure Appl. Anal. 13 (2014), no.~6, 2359-2376.
 
\bibitem{Cheng} M. Cheng,
{\it Bound state for the fractional Schr\"{o}dinger equation with unbounded 
potential.}
J. Math. Phys. 53 (2012), 043507.

\bibitem{Cho} Y. Cho and T. Ozawa,
{\it A note on the existence of a ground state solution to a fractional Schr\"odinger equation}.
Kyushu J. Math. 67 (2013), 227-236.

\bibitem{DdW}J. Davila, M. del Pino and J. Wei,{\it  Concentrating standing waves for the fractional nonlinear Schrodinger equation.}
 J. Diff. Eqns. 256(2014), no.2, 858-892.

\bibitem{davila-delpino-dipierro-valdinoci}
J. D\'avila, M. Del Pino, S. Dipierro and E. Valdinoci,
{\it Concentration phenomena for the nonlocal Schr\"odinger equation with
Dirichlet datum}.
Preprint \\ {\tt http://arxiv.org/pdf/1403.4435.pdf}

\bibitem{guide}
E. Di Nezza, G. Palatucci and E. Valdinoci,
{\it Hitchhiker???s guide to the fractional Sobolev spaces}.
Bull. Sci. Math. 136 (2012), no. 5, 521-573.

\bibitem{dipierro}
S. Dipierro, G. Palatucci and E. Valdinoci,
{\it Existence and symmetry results for a
Schr\"odinger type problem involving
the fractional Laplacian}.
Le Matematiche 68 (2013), no. 1, 201-216.

\bibitem{sDSV}
S. Dipierro, O. Savin and E. Valdinoci,
{\it All functions are locally $s$-harmonic up to a small error}.
Preprint \\ {\tt http://www.wias-berlin.de/preprint/1944/wias\_preprints\_1944.pdf}

\bibitem{GF} G. Ev\'{e}quoz and  M. M. Fall, 
{\it Positive solutions to some nonlinear fractional Schr\"{o}dinger
equation via a Min-Max procedure}. \\
{\tt http://arxiv.org/abs/1312.7068}

\bibitem{FW} M. M. Fall and T. Weth, 
{\it Monotonicity and nonexistence results for some fractional
elliptic problems in the half space}.
To appear in Commun. Contemp. Math.\\ 
{\tt http://arxiv.org/abs/1309.7230}

\bibitem{felmer} P. Felmer, A. Quaas and J. Tan,
{\it Positive solutions of Nonlinear Schr\"{o}dinger 
equation with the fractional Laplacian}. 
Proc. R. Soc. Edinb., Sect. A, Math. 142 (2012), no. 6, 1237-1262.

\bibitem{Feng} B. Feng,
{\it Ground states for the fractional Schr{\"o}dinger equation},
Electron. J. Differential Equations 127 (2013), 1-11.

\bibitem{FLS} R. L. Frank, E. Lenzmann and L. Silvestre,
{\it Uniqueness of radial solutions to fractional laplacian}.
Preprint \\ {\tt http://arxiv.org/abs/1302.2652v1}

\bibitem{JLX} T. Jin, Y. Y. Li and J. Xiong,
{\it On a fractional Nirenberg problem, part I: blow up analysis and
compactness of solutions}. J. Eur. Math. Soc. (JEMS) 16 (2014),
no. 6, 1111-1171.

%\bibitem{Lions} P.-L. Lions,
%{\it Sym\'etrie et compacit\'e dans les espaces de Sobolev}.
%J. Funct. Anal. 49 (1982), no. 3, 315-334.

\bibitem{NiTa93} W-M. Ni; I. Takagi, {\it
Locating the peaks of least-energy solutions to a semilinear Neumann problem.} 
Duke Math. J. 70 (1993), no. 2, 247-281. 

\bibitem{NiTa} W.-M. Ni and I. Takagi, 
{\it On the shape of least-energy solutions to a semilinear Neumann problem}. 
Comm. Pure Appl. Math. 44 (1991), no. 7, 819-851. 

\bibitem{NW} W-M. Ni and J. Wei, {\it On the location and profile of spike-layer solutions
to singularly perturbed semilinear Dirichlet problems}, Comm. Pure Appl. Math. 48, 1995, 731-768.


\bibitem{SW}S. Santra, J. Wei, {\it
Profile of the least energy solution of a singular perturbed Neumann problem with mixed powers. }
Ann. Mat. Pura Appl. (4) 193 (2014), no. 1, 39-70. 
\bibitem{secchi}
S. Simone,
{\it Ground state solutions for nonlinear fractional
Schr\"odinger equations in~$\R^N$}.
J. Math. Phys. 54 (2013), no. 3, 031501, 17 p.

\bibitem{Sil} L. Silvestre,
{\it Regularity of the obstacle problem for a fractional power of the
Laplace operator}. Comm. Pure Appl. Math. 60 (2007), no.~1, 67-112.
   
\bibitem{Wang} X. Wang,
{\it On concentration of positive bound states of nonlinear Schr\"{o}dinger 
equations}, Commun. Math. Phys. 153 (1993), 229-244. \\
\bigskip
 
\end{thebibliography}
\end{document}